\documentclass[final]{siamltex}
 
\usepackage[boxed,linesnumbered]{algorithm2e}
\usepackage{amsmath}
\usepackage{amssymb}
\usepackage{latexsym,pifont,color,comment}
\usepackage{epsfig,graphicx}

\newtheorem{remark}{Remark}
  
\def\A{{A}}
\def\Av{{\bf A}}
\def\B{{\bf B}}

\def\u{{\bf u}}

\def\x{{\bf x}}

\def\f{{\bf f}}
\def\g{{\bf g}}

\def\F{{\bf F}}
\def\G{{\bf G}}

\def\En{{\mathcal E}}

\def\divg{\nabla\cdot}
\def\curl{\nabla \times}

\def\Dx{\Delta x}
\def\Dt{\Delta t}
\def\hf{\frac{1}{2}}

\title{Positivity-Preserving Finite Difference WENO Schemes with Constrained Transport for Ideal Magnetohydrodynamic Equations}

\author{Andrew J. Christlieb\thanks{Department of Mathematics and Department of Electrical and Computer Engineering,
Michigan State University, East Lansing, MI 48824, USA ({\tt christli@msu.edu}).} \and
   Yuan Liu\thanks{Department of Mathematics, Michigan State University, 
Michigan State University, East Lansing, MI 48824, USA ({\tt yliu7@math.msu.edu}).} \and
Qi Tang\thanks{Department of Mathematics, Michigan State University, East Lansing, MI 48824, USA
({\tt tangqi@msu.edu}).} \and
Zhengfu Xu\thanks{Department of Mathematical Science, Michigan Technological University, Houghton, MI 49931, USA
({\tt zhengfux@mtu.edu}).}
}

\begin{document}

\maketitle

\begin{abstract}

In this paper, we utilize the maximum-principle-preserving flux limiting technique, originally designed for high order weighted essentially non-oscillatory (WENO) methods for scalar hyperbolic conservation laws, 
to develop a class of high order positivity-preserving finite difference WENO method for the ideal magnetohydrodynamic (MHD) equations.   Our scheme,  under the constrained transport (CT) framework, can achieve high order  accuracy, a discrete divergence-free condition and positivity of the numerical solution   simultaneously. 
Numerical examples in 1D, 2D and 3D are provided to demonstrate the performance of the proposed method.

\end{abstract}



\section{Introduction}
\label{sec1}

In this paper, we  propose a class of high-order positivity-preserving finite difference WENO schemes within the unstaggered CT framework for the ideal MHD equations. 
The ideal MHD equations are fluid models of perfectly conducting
quasi-neutral plasmas, and consist of nonlinear hyperbolic conservation laws for the macroscopic quantities 
with an additional divergence-free restriction on the magnetic field. 
Mathematically, the ideal MHD equations can be written in a conservative form as follows,
\begin{gather}
\label{MHD}
  \frac{\partial}{\partial t}
   \begin{bmatrix}
     \rho \\ \rho \u \\ \En\\ \B
    \end{bmatrix}
    + \divg{
    	\begin{bmatrix}
	  \rho \u \\ 
	  \rho \u \otimes \u + (p+\frac{1}{2} \| \B \|^2) \mathbb{I} - \B \otimes \B \\
	  \u (\mathcal{E} +  p +\frac{1}{2} \| \B \|^2 ) - \B(\u \cdot \B) \\
	  \u \otimes \B- \B \otimes \u 
	\end{bmatrix}
	}
	 = 0, \\
\label{constraint}
\divg{ \B} = 0,
\end{gather} 
with 
\begin{align}
\label{eq:eng}
	\mathcal{E} = \frac{p}{\gamma-1} + \frac{\rho \|\u\|^2}{2}   + \frac{\|\B\|^2}{2 }.
\end{align}
Here, $\rho$ is density of mass, $\rho \u$ is momentum,  $\mathcal{E}$ is total energy, $p$ is the hydrodynamic pressure, $\| \cdot \|$ is used to denote the Euclidean vector norm and $\gamma = 5/3$ is the ideal gas constant. 

One difficulty to simulate the ideal MHD equations is how to propagate a discrete version of the divergence-free condition forward in time.
The failure to satisfy this condition produces numerical instabilities and has been well documented in the literature \cite{article:BrBa80, article:EvHa88, article:To00}.
To design divergence-free methods for solving the ideal MHD equations, 
the CT methodology arises as one important approach, see \cite{article:Ba04,article:BaSp99a,article:ChRoTa14,article:DaWo98,article:EvHa88,article:FeTo03,article:HeRoTa10,article:LoZa00,article:LoZa04,article:Ro04b,article:Ryu98,article:De01b,article:To00,article:To2005,article:Ro13} for references.
Following  \cite{article:ChRoTa14,article:HeRoTa10,article:helzel2013,article:Ro04b,article:Ro13}, we propose to conduct our investigation within the CT framework in this paper.

Another major focus of this paper is the design of high-order schemes that preserve the positivity 
of the density and pressure of the MHD system.  Even with divergence-free methods, 
negative density or/and pressure can still be observed in   numerical simulations, 
such as those for the low-$\beta$ plasma.  
This can lead to a complex wave speed that breaks the hyperbolicity of the system and 
causes the numerical simulations to break down.  A lot of efforts have been dedicated 
addressing this issue in the literature. For instance, 
Balsara and Spicer \cite{article:BaSp99b} proposed a strategy to maintain the positivity of 
pressure by switching the Riemann solvers based on different wave situations.
Janhunen \cite{article:Ja00} designed a new Riemann solver for the modified ideal MHD equations 
and demonstrated its positivity-preserving property numerically.
In \cite{article:Wa09}, a conservative second-order MUSCL-Hancock
scheme was shown to be positivity-preserving for the 1D ideal MHD equations 
and the extension to multi-dimensional (multi-D) cases was constructed based on  similar ideas as 
Powell's 8-wave formulation \cite{article:Po94,article:Po99}.
Balsara \cite{article:Balsara12}  developed a high-order positivity-preserving scheme for ideal MHD through limiting high-order
numerical solutions by a conservative bounded solution. 
Another class of important methods for the ideal MHD equations is discontinuous Galerkin (DG) methods \cite{article:LiShu05, article:LiXu12, article:LiXuYa11, article:Ro13, article:YaXuLi13}.
Recently, Cheng \emph{et al.} proposed positivity-preserving  DG and central DG methods for the ideal MHD equations \cite{article:Cheng13}, 
in which they generalized Zhang and Shu's positivity-preserving limiters for the compressible Euler equations
\cite{article:ZhangShu10rectmesh}.
In \cite{article:Cheng13},  it was proven that
 the first-order Lax-Fridrichs scheme is positivity-preserving for the 1D MHD under the restriction CFL$\le 0.5$. This first-order scheme also serves as the building block for the positivity-preserving scheme in this paper.


Besides the aforementioned work for MHD equations, several high-order positivity-preserving schemes have been developed recently 
for compressible Euler equations.
Zhang and Shu  developed arbitrary-order positivity-preserving finite volume WENO and DG methods by limiting the underlying  polynomials around cell averages \cite{article:ZhangShu10rectmesh}.
A flux cut-off limiter  was  proposed by Hu \emph{et al.}\ \cite{article:hu2013positivity} for
 finite difference WENO schemes to maintain positivity of density and pressure for the compressible Euler system.
 In this paper, we adopt the parametrized positivity-preserving flux limiter for the compressible Euler systems in \cite{article:XiQiXu}, which was originated from the maximum-principle-preserving flux limiter in \cite{liang2014parametrized,xu2013} for 1D and 2D scalar hyperbolic conservation laws. The approach developed in  \cite{article:ChLiTaXu14,article:XiQiXu, TaQiXu13} is novel because  the parametrized limiter is only applied at the final stage of RK method, making the implementation more efficient and maintaining the accuracy of the base scheme without sacrificing the CFL excessively.
 

The rest of the paper is organized as follows. 
In Section 2 we will briefly review the evolution
equations for the magnetic potential in 2D and 3D MHD systems and the general framework
of the CT approach.
In Section 3 and 4, we present positivity-preserving finite difference WENO schemes for the 1D and multi-D MHD equations. 
The proposed schemes are implemented and tested on several 1D, 2D and 3D numerical 
examples in Section~\ref{numer}. The conclusions are given in Section~\ref{conclu}.

\section{Review of WENO constrained transport schemes}
\label{sec:review}
In this section, we will briefly review the concepts of the magnetic potential in a CT framework and 
 outline the WENO-CT schemes in \cite{article:ChRoTa14}.

\subsection{Magnetic potential}
\label{mhd-eq-mq3d}
In the CT framework, instead of solving the magnetic field directly, a magnetic potential is introduced in order to reconstruct a discrete divergence-free magnetic field. For example, the divergence-free magnetic field can be written as the curl of a magnetic vector potential in the 3D MHD system,
\begin{align}
\label{ptform}
\B = \nabla \times \Av.
\end{align}
Furthermore, because of the relation
\begin{align}
\divg \left( \u \otimes \B - \B \otimes \u\right) = \curl(\B \times \u),
\end{align} 
the magnetic induction equation in $\eqref{MHD}$ can be rewritten 
in curl form:
\begin{align}
\label{bevolve}
\frac{\partial \B}{\partial t}+ \nabla \times (\B \times \u) = 0.
\end{align}
Substituting the magnetic vector potential $\eqref{ptform}$ into the evolution equation $\eqref{bevolve}$, we  obtain
\begin{align}
\label{curlA}
\curl \left\{ \frac{\partial \Av}{\partial t}+ (\curl \Av)  \times \u \right\} = 0.
\end{align}
Therefore, there exists a scalar potential function $\psi$ such that
\begin{align}
\frac{\partial \Av}{\partial t} + (\curl \Av)  \times \u = -\nabla \psi.
\end{align}
An extra gauge condition is needed to uniquely determine the potential function $\psi$.

Helzel et al.\ \cite{article:HeRoTa10}  investigated different choices of gauge conditions and found that  stable solutions can be obtained by introducing the Weyl gauge, i.e., setting $\psi \equiv 0$. With this gauge condition, the evolution equation for the vector potential becomes
\begin{align}
\label{3dpte}
\frac{\partial \Av}{\partial t} + (\curl \Av)  \times \u = 0.
\end{align}
We notice that the 2D MHD system actually results in a simpler version of \eqref{3dpte}, because the divergence-free condition is reduced to
\begin{align}
\label{2dcons}
\nabla \cdot \B = \frac{\partial B_x}{\partial x} + \frac{\partial  B_y}{\partial y} = 0,
\end{align}
where $B_x$ and $B_y$ are reconstructed with only the third component of the magnetic potential,
\begin{align}
\label {2dcurl}
B_x =  \frac{\partial A_z}{\partial y} \quad \text{and} \quad B_y = - \frac{\partial A_z}{\partial x},
\end{align}
 effectively reducing the vector potential $\Av$ to a {scalar} potential $A_z$. In this case, \eqref{3dpte} is  reduced to
\begin{align}
\label{2dpte}
 \frac{\partial A_z}{\partial t} + u_x  \frac{\partial A_z}{\partial x} + u_y \frac{\partial A_z}{\partial y} = 0.
\end{align}

It is worthwhile to point out that the full vector potential evolution equation $\eqref{3dpte}$ is a non-conservative, weakly-hyperbolic system while the scalar potential equation $\eqref{2dpte}$ is strongly hyperbolic. In \cite{article:ChRoTa14}, Christlieb \emph{et al.}\ proposed a class of finite difference schemes based on WENO reconstruction to solve both the scalar potential in 2D and vector potential in 3D. In particular, the authors introduced   an artificial resistivity approach for the 3D system $\eqref{3dpte}$ in order to control the unphysical oscillations in the magnetic field.

\subsection{Outline of WENO constrained transport schemes}
\label{sec:tang}
In this subsection, we will present an 
outline of the fundamental CT framework detailed as follows. 
 \medskip
 
A single time-step of the WENO-CT method from time $t^n$ to time $t^{n+1}$ consists of the following sub-steps:

\begin{enumerate}
\setcounter{enumi}{-1}
\item Start with $(\rho^{n}, \rho \u^{n}, \En^{n}, \B^{n})$  and $\text{A}^n$, where $\text{A}^n$ stands for $\Av$ in 3D and $\A_z$ in 2D at time $t^n$.

\item Discretize the MHD equations $\eqref{MHD}$ for the conserved quantities and the potential equation $\eqref{3dpte}$ or $\eqref{2dpte}$ for the magnetic potential by using finite difference WENO schemes in \cite{article:ChRoTa14} and the  strong stability-preserving Runge-Kutta  (SSP-RK) time-stepping method \cite{gottliebShuTadmor01}. This  updates the conserved quantities and the magnetic potential by
\begin{align}
\label{hcl:step2}
(\rho^{n}, \rho \u^{n}, \En^{n}, \B^{n}) & \Rightarrow (\rho^{n+1}, \rho \u^{n+1}, \En^{*}, \B^{*}), \\
\text{A}^{n} & \Rightarrow \text{A}^{n+1},
\end{align}
where $\B^{*}$ is the predicted magnetic field that is not necessarily discrete divergence-free and $\En^{*}$ is the predicted energy.
\item Correct $\B^{*}$ by computing a discrete curl of the magnetic potential $\text{A}^{n+1}$:
\begin{align}
\B^{n+1} = \curl \text{A}^{n+1}.
\end{align}
\item Set the corrected total energy density $\En^{n+1}$ based on one of the following options:
	{\begin{enumerate}
	\item[] \textbf{Option 1:} Conserve the total energy: 
								\begin{align}
								\En^{n+1} = \En^{*}.
								\end{align}
	\item[] \textbf{Option 2:} Keep the pressure the same before and after the magnetic field 
		correction step ($p^{n+1} = p^{*}$): 
								\begin{align}
								\En^{n+1} = \En^{*}+\frac{1}{2} \left( \|\B^{n+1}\|^2 -\|\B^{*}\|^2 \right).
								\end{align}
	\end{enumerate}}
\end{enumerate}
Depending on the scalar or vector magnetic potential used, we call the overall scheme as WENO-CT2D or WENO-CT3D. 

In this paper, we make exclusive use of \textbf{Option 2} in order to preserve the positivity of the pressure after the magnetic field is corrected 
albeit at the expense of sacrificing the energy conservation. 
This is a common technique in the CT framework
 for problems  involving very low $\beta$ plasma  \cite{article:BaSp99a, article:To00}. 
Under this option, if the density and pressure after Step 1 are non-negative, they will be non-negative in the overall computation. Therefore, in numerical computations, it suffices to restrict our attention to designing    positivity-preserving schemes for  \eqref{hcl:step2} in Step 1.

Another difference of the schemes considered in this paper compared to those in \cite{article:ChRoTa14} lies in the implementation of the correction steps (Steps 2 and 3). We propose to perform the correction steps only at the end of each time step $t^n$   instead of each stage of RK methods in \cite{article:ChRoTa14}.
By doing this modification, we can focus on the final stage of the solution when implementing the limiting technique in Section~\ref{sec:1d}.
Numerical results show negligible differences between the two approaches when SSP-RK3 time-stepping is used.
However, we note that this modification may result in accumulation of the divergence error especially for RK methods with large stage numbers, such as the low-storage 10-stage SSP-RK4 method considered in \cite{article:ChRoTa14}. For those time stepping schemes, this kind of modification is not recommended and the correction steps have to be performed at each stage.


\section{1D case}
\label{sec:1d}
In this section, we describe our positivity-preserving scheme on a 1D MHD system.
The divergence-free condition $\divg \B = 0$ in 1D case is equivalent to $B_x = \text{constant}$.
Since the WENO hyperbolic conservation law solver (WENO-HCL) in \cite{article:JiShu96, article:JiWu99} without CT approaches will produce a solution with constant $B_x$, we use it as our MHD base scheme in 1D, to which we apply a positivity-preserving limiter.

The MHD equations $\eqref{MHD}$ in 1D can be written as follows:
\begin{align}
\label{eqn:conslaw}
\frac{\partial q}{\partial t} + \frac{\partial}{\partial x} \f (q) = 0,
\end{align}
where
\begin{gather}
q =  \left(\rho, \rho u_x, \rho u_y, \rho u_z, \mathcal{E}, B_x, B_y,B_z\right), \\
\begin{split}
\f(q)  &=  \biggl( \rho u_x, \rho u_x u_x + p + \frac{1}{2} \|\B\|^2 - B_x B_x, \rho u_x u_y -  B_x B_y,  \rho u_x u_z - B_x B_z, \\ 
  & \qquad u_x \left(\mathcal{E}+p+\frac{1}{2}\|\B\|^2 \right)- B_x(\u \cdot \B), 0, u_x B_y - u_y B_x, u_x B_z - u_z B_x \biggr).
\end{split}
\end{gather}

The spatial domain $[0,1]$ is divided into $N$ uniform cells:
\begin{equation}\label{1d_eq2}
0 = x_{\hf} < x_{\frac{3}{2}} < ...< x_{N+\hf} =1,
\end{equation}
and we denote
\begin{align*}
I_j = [x_{j-\hf}, x_{j+\hf}], \qquad x_j = \hf( x_{j-\hf}+ x_{j+\hf}), \qquad
\Dx_j=\Dx = 1/N.
\end{align*}
Let ${q}_j(t)$ be the numerical solution at the grid point $x_j = \hf( x_{j-\hf}+ x_{j+\hf}) $. 
The finite difference WENO-HCL schemes solve  $\eqref{eqn:conslaw}$ by a conservative form:
\begin{equation}\label{1d_eq3}
\frac{d}{dt} {q}_j(t) +\frac{1}{\Delta x}\left(\hat{\F}_{j+\hf}-\hat{\F}_{j-\hf} \right) = 0,
\end{equation}
where $\hat{\F}_{j+\hf}$ is defined as a high-order numerical flux constructed by WENO-HCL.
The design of $\hat{\F}_{j+\hf}$ involves
eigenvalue decompositions, physical flux splitting and WENO reconstruction,
the details of which can be found in many references such as \cite{article:ChRoTa14, article:JiWu99}.
One numerical difficulty is that the wave speeds of the MHD system involve the term $1/\rho$.
To avoid the possibility of an infinite wave speed during the computation,
we assume there is a small lower bound $\epsilon_0$ for both density and pressure in the exact solution of the problem we considered.

The semi-discrete equation $\eqref{1d_eq3}$ can be further discretized in time by high-order time integrators. While  our proposed scheme can be applied with any RK method,  
we take the following third-order SSP-RK method as an illustrative example:
\begin{align}
\label{tvdrk}
\begin{split}
{q}_j^{(1)} & = {q}_j^n + \Dt L({q}_j^n), \\
{q}_j^{(2)} & = {q}_j^n + \frac{1}{4}\Dt \left(L({q}_j^n)+L(q_j^{(1)})\right), \\
{q}_j^{n+1} & = {q}_j^n + \frac{1}{6}\Dt \left(L(q_j^n)+4L(q_j^{(2)})+L(q_j^{(1)})\right),
\end{split}
\end{align}
where ${q}_j^{(k)}$ and ${q}_j^n$ denote the numerical solutions at the $k^\text{th}$ RK stage and $t=t^n$ respectively, and
\begin{align}
L(q^n_j) = -\frac{1}{\Dx}(\hat{\F}^n_{j+\hf} -\hat{\F}^n_{j-\hf}).
\end{align}
If we use $\hat{\F}^{n}_{j+\hf}$, $\hat{\F}^{(1)}_{j+\hf}$ and $\hat{\F}^{(2)}_{j+\hf}$
to denote the numerical fluxes reconstructed based on $q^n$, $q^{(1)}$ and $q^{(2)}$,
the final stage of RK discretization $\eqref{tvdrk}$ can be rewritten as,
\begin{equation}\label{1d_eq4}
q^{n+1}_j = q^n_j - \lambda(\hat{\F}^{rk}_{j+\hf}-\hat{\F}^{rk}_{j-\hf}) ,
\end{equation}
where
\begin{align}
\lambda = \frac{\Dt}{\Dx}, \qquad \hat{\F}^{rk}_{j+\hf} =\frac{1}{6}\left( \hat{\F}^n_{j+\hf} + 4\hat{\F}^{(2)}_{j+\hf} +\hat{\F}^{(1)}_{j+\hf} \right).
\end{align}
$\hat{\F}^{rk}_{j+\hf}$ can be viewed as a linear combination of high-order numerical fluxes from different stages. Following the ideas in \cite{article:XiQiXu,xu2013}, we need to modify the numerical flux  $\hat{\F}^{rk}_{j+\hf}$ by a positivity-preserving flux to design a high-order positivity-preserving MHD scheme.

Cheng \emph{et al.}\ \cite{article:Cheng13} proved the simple Lax-Friedrichs numerical flux coupled with forward Euler time discretization
is positivity-preserving for the 1D MHD equations $\eqref{eqn:conslaw}$ under the restriction
CFL$\le 0.5$. 
When the Lax-Friedrichs scheme is used to solve the high-order solution $q^n$ from $t^n$ to $t^{n+1}$,
we have
\begin{align}
\label{liu1}
\hat{q}^{n+1}_j = {q}^n_j - \lambda(\hat{\f}_{j+\hf}-\hat{\f}_{j-\hf}),
\end{align}
where $\hat{q}^{n+1}_j$ is introduced to denote the low-order solution at $x_j$ and $t = t^{n+1}$, and the Lax-Friedrichs flux is formulated as
\begin{align}
\label{eqn:lf}
\hat{\f}_{j+\hf} = \hf \left( \f({q}_{j+1}^n)+ \f({q}_j^n) - \alpha ({q}_{j+1}^n-{q}_j^n) \right),
\end{align}
where the maximal wave speed $\alpha$ is defined by,
\begin{align}
\alpha = \max_x \left( |u_x|+c_x^f \right).
\end{align}
Here $c_x^f$ is the fast speed of the MHD system in x-direction, see \cite{article:Po94} for reference.

The density and pressure  computed by the first-order scheme \eqref{liu1} satisfy  
\begin{align}
\label{eqn:pp1}
\begin{cases}
\hat{\rho}^{n+1}_j > 0, \\
\hat{p}^{n+1}_j >0.
\end{cases}
\end{align}
We can use the first-order solution $\hat{q}^{n+1}$ to define 
the numerical lower bounds for the density and pressure for the high-order solution ${q}^{n+1}$, 
which are 
\begin{align}
\epsilon_\rho^{n+1} & = \min_j (\hat{\rho}_j^{n+1},\epsilon_0), \label{eps_rho} \\
\epsilon_p^{n+1} & = \min_j (\hat{p}_j^{n+1},\epsilon_0). \label{eps_p}
\end{align}
Throughout the simulations for this work, we take $\epsilon_0 = 10^{-13}$.
It can be certainly taken as a smaller number if it is required by the problem and allowed by the machine precision.

Following \cite{article:XiQiXu,xu2013}, to guarantee the positivity of the  high-order solutions by the WENO scheme $\eqref{1d_eq4}$, we need to find a modification of the numerical flux as follows:
\begin{align}
\tilde{\F}_{j+\hf} = \theta_{j+\hf}(\hat{\F}^{rk}_{j+\hf} - \hat{\f}_{j+\hf}) + \hat{\f}_{j+\hf}.
\end{align}
where the limiting parameter $\theta_{j+\hf} \in [0,1]$. 
We seek a combination of $\theta_{j+\hf}$, such that the solutions satisfy
\begin{align}
\label{ppsimu}
\begin{cases}
{\rho}^{n+1}_j \geq \epsilon_\rho^{n+1}, \\
{p}^{n+1}_j \geq \epsilon_p^{n+1}. 
\end{cases}
\end{align}

Our positivity-preserving limiting technique follows a two-step procedure. Firstly, as outlined below, we describe a strategy to guarantee the computed density positive. To facilitate the discussion, we denote the first order flux of the density in $\hat{\f}$ as  $\hat{f}^\rho$, whereas ${f}^\rho$ and $\tilde{f}^\rho$ are the corresponding flux components in $\hat{\F}^{rk}$ and $\tilde{\F}$, respectively.

To preserve positive density, we need to find  upper bounds $\Lambda^\rho_{\pm\hf, I_j}$ of the limiting parameters $\theta_{j \pm \hf}$ at each cell $I_j$, such that, for any combination
$ (\theta_{j-\hf}, \theta_{j+\hf}) \in [0, \Lambda^\rho_{-\hf, I_j}] \times [0, \Lambda^\rho_{+\hf, I_j}]$,
the following inequality holds:
\begin{align}
\label{liu3}
{\rho}_j^{n+1} (\theta_{j-\hf},\theta_{j+\hf}) = {\rho}_j^n - \lambda ( \tilde{f}^\rho_{j+\hf} - \tilde{f}^\rho_{j-\hf} ) \geq \epsilon_\rho^{n+1},
\end{align}
where $\tilde{f}^\rho_{j+\hf} = \theta_{j+\hf}(f^\rho_{j+\hf} - \hat{f}^\rho_{j+\hf}) + \hat{f}^\rho_{j+\hf},$
Further introducing the notation $\Gamma_j = {\rho}_j^n - \lambda ( \hat{f}^\rho_{j+\hf} - \hat{f}^\rho_{j-\hf} )$,
\eqref{liu3} is equivalent to
\begin{align}
\label{inequ1}
\Gamma_j -  \lambda( \theta_{j+\hf}({f}^\rho_{j+\hf} - \hat{f}^\rho_{j+\hf}) 
- \theta_{j-\hf}({f}^\rho_{j-\hf} - \hat{f}^\rho_{j-\hf})  ) \geq \epsilon_\rho^{n+1}.
\end{align}
Due to the positivity-preserving property of the first-order scheme and the definition of $\epsilon_\rho^{n+1}$ $\eqref{eps_rho}$, 
we have $\Gamma_j \ge \epsilon_\rho^{n+1}$.
Thus, the inequality $\eqref{inequ1}$ can be rewritten as,
\begin{align}
\label{inequ2}
 \lambda \theta_{j-\hf}({f}^\rho_{j-\hf} - \hat{f}^\rho_{j-\hf}) 
- \lambda \theta_{j+\hf}({f}^\rho_{j+\hf} - \hat{f}^\rho_{j+\hf})   \geq \epsilon_\rho^{n+1} - \Gamma_j
\end{align}
with the right hand side $\epsilon_\rho^{n+1} - \Gamma_j \le 0$. For abbreviation, we introduce  a notation $F_{j+\hf} = {f}^\rho_{j+\hf} - \hat{f}^\rho_{j+\hf}$.

Following the same idea in \cite{article:XiQiXu,xu2013}, we will determine
the upper bounds of the parameter $\theta_{j\pm\hf}$ by a case-by-case discussion 
based on the signs of $F_{j-\hf}$ and $F_{j+\hf}$. 
In particular, we decouple the inequalities $\eqref{inequ2}$ based on the following four cases:
\begin{itemize}
  \item If $F_{j-\hf} \ge 0$ and $F_{j+\hf} \le 0$, then
    $$(\Lambda^\rho_{-\hf, I_j}, \Lambda^\rho_{+\hf, I_j}) = (1,1). $$
  \item If $F_{j-\hf} \ge 0$ and $F_{j+\hf} > 0$, then
    $$(\Lambda^\rho_{-\hf, I_j}, \Lambda^\rho_{+\hf, I_j}) = \left(1,\min \left(1, \frac{\epsilon_\rho^{n+1}-
    \Gamma_j}{-\lambda F_{j+\hf}}\right)\right). $$
  \item If $F_{j-\hf} < 0$ and $F_{j+\hf} \le 0$, then
  $$(\Lambda^\rho_{-\hf, I_j}, \Lambda^\rho_{+\hf, I_j}) = \left(\min \left(1, \frac{\epsilon_\rho^{n+1}-\Gamma_j}{\lambda F_{j-\hf}}\right),1 \right). $$

  \item If $F_{j-\hf} < 0$ and $F_{j+\hf} >0$,
  \begin{itemize}
    \item  if the inequality $\eqref{inequ2}$ is satisfied with $(\theta_{j-\hf}, \theta_{j+\hf}) = (1, 1)$ then
    $$(\Lambda^\rho_{-\hf, I_j}, \Lambda^\rho_{+\hf, I_j}) = (1,1). $$
    \item  otherwise, we choose
    $$(\Lambda^\rho_{-\hf, I_j}, \Lambda^\rho_{+\hf, I_j}) = \left(\frac{\epsilon_\rho^{n+1} - \Gamma_j}{\lambda F_{j-\hf} - \lambda F_{j+\hf}},\frac{\epsilon_\rho^{n+1} - \Gamma_j}{\lambda F_{j-\hf} - \lambda F_{j+\hf}} \right). $$
  \end{itemize}
  \end{itemize}
This procedure has been discussed in \cite{article:XiQiXu,xu2013}.
It is easy to show when $ (\theta_{j-\hf}, \theta_{j+\hf}) \in [0, \Lambda^\rho_{-\hf, I_j}] \times [0, \Lambda^\rho_{+\hf, I_j}]$ with the bounds $\Lambda^\rho_{\pm\hf, I_j}$ obtained by the above strategy, the inequality $\eqref{inequ2}$ holds, i.e., the density $\rho^{n+1}_j$ is positive at each grid $x_j$. We define this set as $S_{\rho,I_j}$:
\begin{align}
S_{\rho,I_j} = [0, \Lambda^\rho_{-\hf, I_j}] \times [0, \Lambda^\rho_{+\hf, I_j}].
\end{align}

We next describe a strategy to obtain positive pressure. 
First we discuss some properties of the pressure function
\begin{align}
p(q) = (\gamma - 1)
\left( {\En} -\hf \frac{(\rho u_x)^2 + (\rho u_y)^2+ (\rho u_z)^2}{\rho}  - \hf \left( B_x^2+B_y^2+B_z^2 \right) \right).
\end{align}
We note the pressure function is concave with respect to $q = \left(\rho, \rho u_x, \rho u_y, \rho u_z, \mathcal{E}, B_x, B_y,B_z\right)$.

Similar as the function $\rho_j^{n+1}(\theta_{j-\hf}, \theta_{j+\hf})$, we can define a function $p_j^{n+1}(\theta_{j-\hf}, \theta_{j+\hf})$ as follows,
\begin{align}
p_j^{n+1}(\theta_{j-\hf}, \theta_{j+\hf}) : = p(q_j^{n+1}(\theta_{j-\hf}, \theta_{j+\hf}))
\end{align}
We need the following lemma to construct the limiter.
\begin{lemma}
\label{lemma1}
The pressure function satisfies
\begin{align}
    p\left( q^{n+1}_j \left( \alpha \overrightarrow{\theta}^1 + (1-\alpha) \overrightarrow{\theta}^2 \right) \right)
    \geq \alpha   p\left( q^{n+1}_j \left( \overrightarrow{\theta}^1 \right) \right)
    + (1- \alpha) p\left( q^{n+1}_j \left( \overrightarrow{\theta}^2 \right) \right)
\end{align}
for any $\alpha \in [0,1]$ and $\overrightarrow{\theta}^1, \overrightarrow{\theta}^2 \in S_{\rho,I_j}$.
\end{lemma}

The proof of this lemma is straightforward, as long as we use the concave property of $p(q)$ and note that the solution $q^{n+1}_j$ is a linear function of its limiting parameters, i.e.,
\begin{align*}
    q^{n+1}_j       \left( \alpha \overrightarrow{\theta}^1 + (1-\alpha) \overrightarrow{\theta}^2 \right) =
    \alpha q^{n+1}_j\left( \overrightarrow{\theta}^1\right) + (1-\alpha) q^{n+1}_j\left(\overrightarrow{\theta}^2 \right).
\end{align*}
A similar lemma in the Euler equations has been use in the past \cite{article:ChLiTaXu14,article:XiQiXu}.

We want to identify a subset of the set $S_{\rho,I_j}$, denoted by $S_{p,I_j}$,
such that $p^{n+1}_j(\theta_{j-\hf},\theta_{j+\hf})$ is positive, i.e.,
\begin{align}
S_{p,I_j} = \{ (\theta_{j-\hf},\theta_{j+\hf} ) \in [0, \Lambda^\rho_{-\hf, I_j}] \times [0, \Lambda^\rho_{+\hf, I_j}]: p^{n+1}_j(\theta_{j-\hf},\theta_{j+\hf}) \ge \epsilon_p^{n+1} \}.
\end{align}
Due to Lemma~\ref{lemma1}, $S_{p,I_j}$ a convex set. To determine $S_{p,I_j}$, we can only focus on its vertices.

If we denote the four vertices of $S_{\rho,I_j}$ to be 
$
A^{k_1,k_2} = (k_1 \Lambda^\rho_{-\hf,j}, k_2 \Lambda^\rho_{+\hf,j}),
$
with $k_1$, $k_2$ being 0 or 1, similarly we can define   the  vertices of $S_{p,I_j}$ to be 
$
B^{k_1,k_2}.
$
For $(k_1,k_2) \neq (0,0)$, if $p_j^{n+1}(A^{k_1,k_2}) \geq \epsilon_p^{n+1}$, we let $B^{{k_1,k_2}} = A^{k_1,k_2}$; otherwise we find $r$ such that $p_j^{n+1}(rA^{k_1,k_2}) \geq \epsilon_p^{n+1}$ and let $B^{{k_1,k_2}} = r A^{k_1,k_2}$. The resulting three vertices $B^{k_1,k_2}$ with the origin $(0,0)$ form $S_{p,I_j}$.

Next, we can identify  a rectangle inside $S_{p,I_j}$ denoted by
\begin{align}
R_{\rho,p,I_j} = [0, \Lambda_{-\hf,I_j}] \times [0,\Lambda_{+\hf,I_j}],
\end{align}
where
\begin{align}
\Lambda_{-\hf,I_j} = \min_{\substack{k_2 = 0, 1 }}(B^{1,k_2}), \quad
\Lambda_{+\hf,I_j} = \min_{\substack{k_1 = 0, 1 }}(B^{k_1,1}).
\end{align}

After repeating this procedure for all $j$, we let
\begin{align}
\theta_{j+\hf} = \min(\Lambda_{+\hf,I_j},\Lambda_{-\hf, I_{j+1}}),
\end{align}
and this finishes our discussion for the 1D MHD scheme.


\begin{remark}
The limiting technique here is only used to guarantee the positivity of the solution at the final stage of Runge-Kutta methods. If there is negative density or pressure in the intermediate stage, we take the absolute value of the density and pressure in the code where a positive solution is required.
The first place that needs a positive solution is to estimate the speed waves of the system. For instance, the speed of sound is taken as $c = \sqrt{\gamma {|p|}/{|\rho|}}$ in the intermediate stage. The second place requiring a positive solution is to estimate the eigenvectors of the Jacobian matrix of the flux function. Those treatments will not degrade the order of accuracy, because the WENO algorithm only needs an estimate of the local eigenvalues and eigenvectors and we always use the true solution to compute the numerical flux even when it becomes negative in the intermediate stage.
However, we also remark that the limiter can be applied to each stage when the positivity in the intermediate stage is required.
\end{remark}

\begin{remark}
From the limiting steps, we can see 
the overall scheme have a CFL constraint of 0.5, which is same as the Lax-Friedrichs scheme.
When the above limiting technique is applied to the intermediate stages, there is no extra restriction because
the time step of the intermediate stage is typically no greater than $\Delta t$.
\end{remark}

\begin{remark}
One numerical difficulty is to satisfy $p(r A^{k_1,k_2}) \ge \epsilon^{n+1}_p$.
This can be done by solving a root $r$ for the equation $p(r A^{k_1,k_2}) = \epsilon^{n+1}_p$.
Through a simple derivation, it can be easily shown the solution $q_j^{n+1}(r A^{k_1,k_2})$ satisfies,
\begin{align*}
q^{n+1}_j(r A^{k_1,k_2}) = r q^{n+1}_j(A^{k_1,k_2}) + (1-r) \hat{q}^{n+1}_j,
\end{align*}
where $\hat{q}^{n+1}_j$ is again used to denote the solution solved by the first order flux $\hat{\f}_{j+\hf}$.
This property is independent from the dimension, which makes it naturally extendable for the multi-D cases.
More importantly, $q^{n+1}_j(A^{k_1,k_2})$ and $\hat{q}^{n+1}_j$ are both computationally cheap to evaluated. 
So with $q^{n+1}_j(A^{k_1,k_2})$ and $\hat{q}^{n+1}_j$ known,
we can solve a root $r$ for the equation,
\begin{align*}
p(r q^{n+1}_j(A^{k_1,k_2}) + (1-r) \hat{q}^{n+1}_j) = \epsilon^{n+1}_p.
\end{align*}
In the MHD equation case, this equation is a cubic function of $r$ in general.
We note that there exist at least one root in the interval $[0, 1]$, 
which can always be found by Newton iteration. 
However, in the implementation,
we only used a simple bisection method with a maximum of ten iterations to find the root, 
because our purpose is to obtain a positive pressure $p(r q^{n+1}_j(A^{k_1,k_2}))$ instead of finding an accurate $r$. 
During the numerical simulations, 
we found the effect of number of iterations to the solution quality and accuracy is negligible. 
A similar approach to find a limiting parameter in positivity-preserving MHD schemes can be found in \cite{article:Balsara12}.
\end{remark}


\section{Multi-D case}
\label{sec:md}
In this section, we briefly describe our positivity-preserving scheme in the multi-D case.
To control the divergence error,  our base scheme is taken as the WENO-CT scheme proposed in \cite{article:ChRoTa14} and outlined in Section \ref{sec:tang}. In the discussion below, we only present the scheme for 2D MHD systems, keeping in mind that the extension to 3D case is quite straightforward. 

The 2D MHD system $\eqref{MHD}$ can be rewritten as:
\begin{align}
\label{eqn:conslaw2}
\frac{\partial q}{\partial t} + \frac{\partial}{\partial x}\f(q) + \frac{\partial}{\partial y} \g(q) = 0.
\end{align}
We need to solve  \eqref{eqn:conslaw2} to get the update in \eqref{hcl:step2}.
If the SSP-RK3 method is used as the time integrator, the WENO-HCL scheme solve
the equation $\eqref{eqn:conslaw2}$ by a conservative form:
\begin{align}
q^{n+1}_{i,j} = q^n_{i,j} - \lambda_x(\hat{\F}^{rk}_{i+\hf,j}-\hat{\F}^{rk}_{i-\hf,j}) 
- \lambda_y(\hat{\G}^{rk}_{i,j+\hf}-\hat{\G}^{rk}_{i,j-\hf}),
\end{align}
where $\hat{\F}^{rk}$ and $\hat{\G}^{rk}$ are linear combinations of high-order numerical fluxes from three RK stages.
Let $\f_{i+\hf,j}$ and $\g_{i,j+\hf}$ again be the first-order Lax-Friedrichs fluxes.
Then we modify the high-order numerical fluxes $\hat{\F}^{rk}$ and $\hat{\G}^{rk}$ by
the Lax-Friedrichs fluxes $\hat{\f}_{i+\hf,j}$ and $\hat{\g}_{i,j+\hf}$ to achieve the positivity of 
the solution, i.e.\
\begin{align}
\tilde{\F}_{i+\hf,j} & = \theta_{i+\hf,j}(\hat{\F}^{rk}_{i+\hf,j} - \hat{\f}_{i+\hf,j}) + \hat{\f}_{i+\hf,j}, \\
\tilde{\G}_{i,j+\hf} & = \theta_{i,j+\hf}(\hat{\G}^{rk}_{i,j+\hf} - \hat{\g}_{i,j+\hf}) + \hat{\g}_{i,j+\hf}.
\end{align}

For each grid $x_{i,j}$, following a two-step strategy similar as the 1D case, we can find a rectangular set
$R_{\rho,p,I_{i,j}} = [0, \Lambda_{L,I_{i,j}}] \times [0,\Lambda_{R,I_{i,j}}] \times [0, \Lambda_{D,I_{i,j}}] \times [0,\Lambda_{U,I_{i,j}}]$, 
such that for any $ (\theta_{i-\hf,j},\theta_{i+\hf,j},\theta_{i,j-\hf},\theta_{i,j+\hf}) \in R_{\rho,p,I_{i,j}}$, we have,
\begin{align}
\rho^{n+1}_{i,j} (\theta_{i-\hf,j},\theta_{i+\hf,j},\theta_{i,j-\hf},\theta_{i,j+\hf}) \ge \epsilon^{n+1}_\rho, \\
p^{n+1}_{i,j}(\theta_{i-\hf,j},\theta_{i+\hf,j},\theta_{i,j-\hf},\theta_{i,j+\hf}) \ge \epsilon^{n+1}_p.
\end{align}
Here $\epsilon^{n+1}_\rho$ and $\epsilon^{n+1}_p$ are the 2D lower bounds with similar definitions as the 1D case 
$\eqref{eps_rho}$ and $\eqref{eps_p}$.
The strategy to find the set $R_{\rho,p,I_{i,j}}$ is similar to the Euler equations case \cite{article:XiQiXu}.
We omit the details here.
After repeating this procedure for all nodes $(i, j)$, we let
\begin{align}
\theta_{i+\hf,j} = \min(\Lambda_{R,I_{i,j}},\Lambda_{L, I_{i+1,j}}), \\
\theta_{i,j+\hf} = \min(\Lambda_{U,I_{i,j}},\Lambda_{D, I_{i,j+1}}).
\end{align}

This whole procedure will produce numerical solution with positive density and pressure after Step 1 in CT framework.
Followed by Step 2 and 3 with {\bf Option 2}, we achieve high order  accuracy, a discrete divergence-free condition and positivity of the numerical solution simultaneously. 
The overall scheme shares the same CFL constraint as the low-order Lax-Fridrichs scheme.
There is no extra restrictions from the limiting process.

As pointed out in \cite{article:Cheng13}, there is still no rigorous proof that the Lax-Friedrichs scheme or any other 
first-order scheme is positivity-preserving in the mutli-D case when the divergence-free constraint is considered.
In this work, we still use the first-order Lax-Friedrichs scheme as the low-order correction scheme for the multi-D cases.
Same as \cite{article:Cheng13}, we take $\text{CFL} \le 0.5$ as the constraint for the positivity-preserving Lax-Friedrichs scheme in the multi-D cases.
On the other hand, our limiting technique is independent from the choice of the low-order scheme.
The overall scheme will be improved as long as we find a positivity-preserving scheme as the building block.

\section{Numerical examples}
\label{numer}
In this section, we perform numerical simulations with our positivity-preserving schemes in 1D, 2D and 3D. SSP-RK3 scheme serves as the time integrator in all the examples whereas fifth-order finite difference WENO-HCL scheme is used for solving the base MHD equations in different examples.
In multi-D, a fourth-order CT method is used to obtain a divergence-free magnetic field.
Unless otherwise stated, the gas constant is $\gamma = 5/3$ and the CFL number is $0.5$.

\subsection{Test cases in 1D}
In this subsection, we test our positivity-preserving scheme by several 1D MHD examples. 
We note that, for all the cases presented in this subsection, negative pressure or density is observed if the base MHD scheme is applied without a positivity-preserving limiter.
Here, the base MHD scheme is fifth-order WENO-HCL scheme.

\subsubsection{Vacuum shock tube test}
\begin{figure}
\begin{center}
\begin{tabular}{cc}
  	(a)\includegraphics[width=0.42\textwidth]{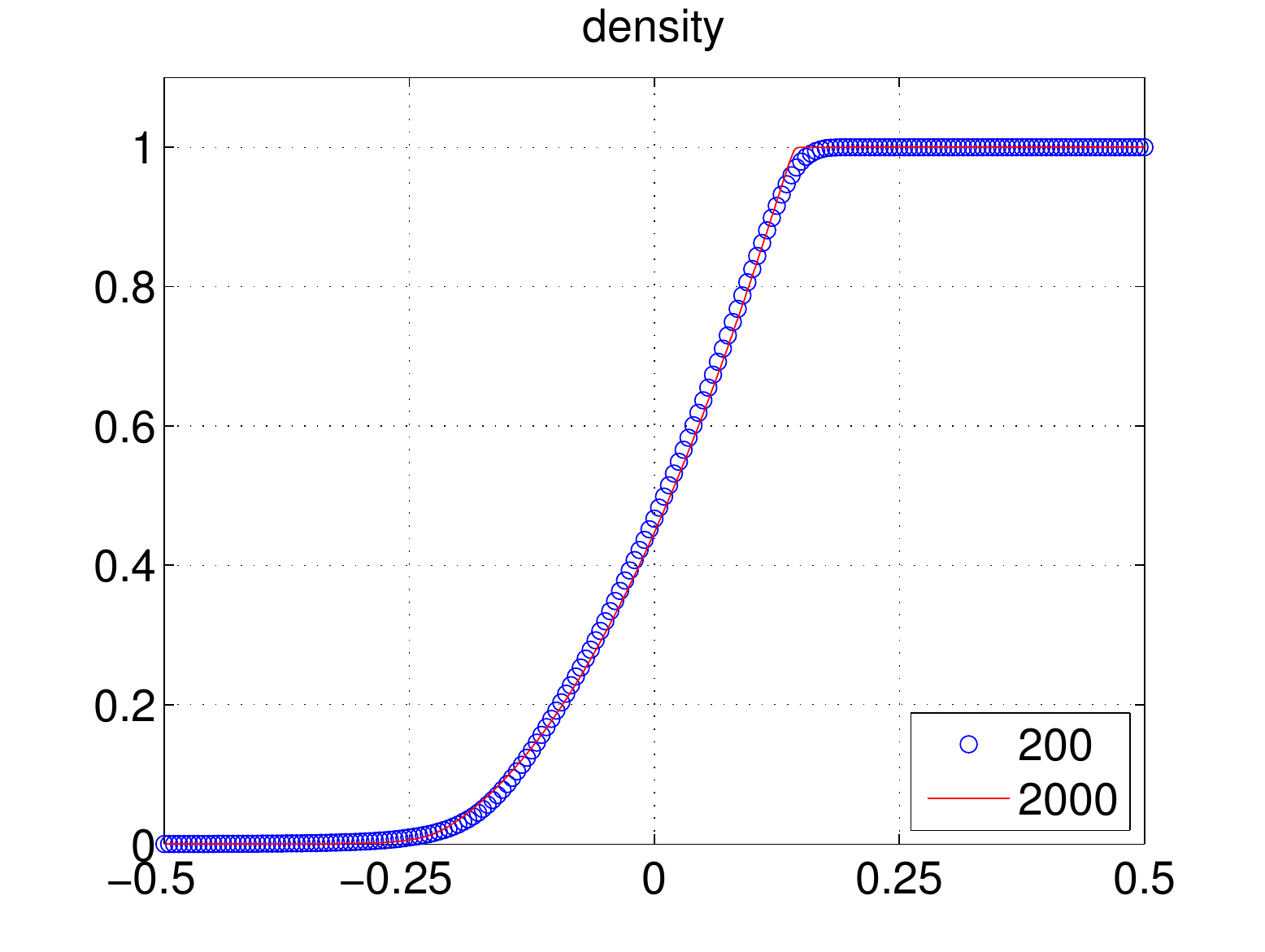} & 
	(b)\includegraphics[width=0.42\textwidth]{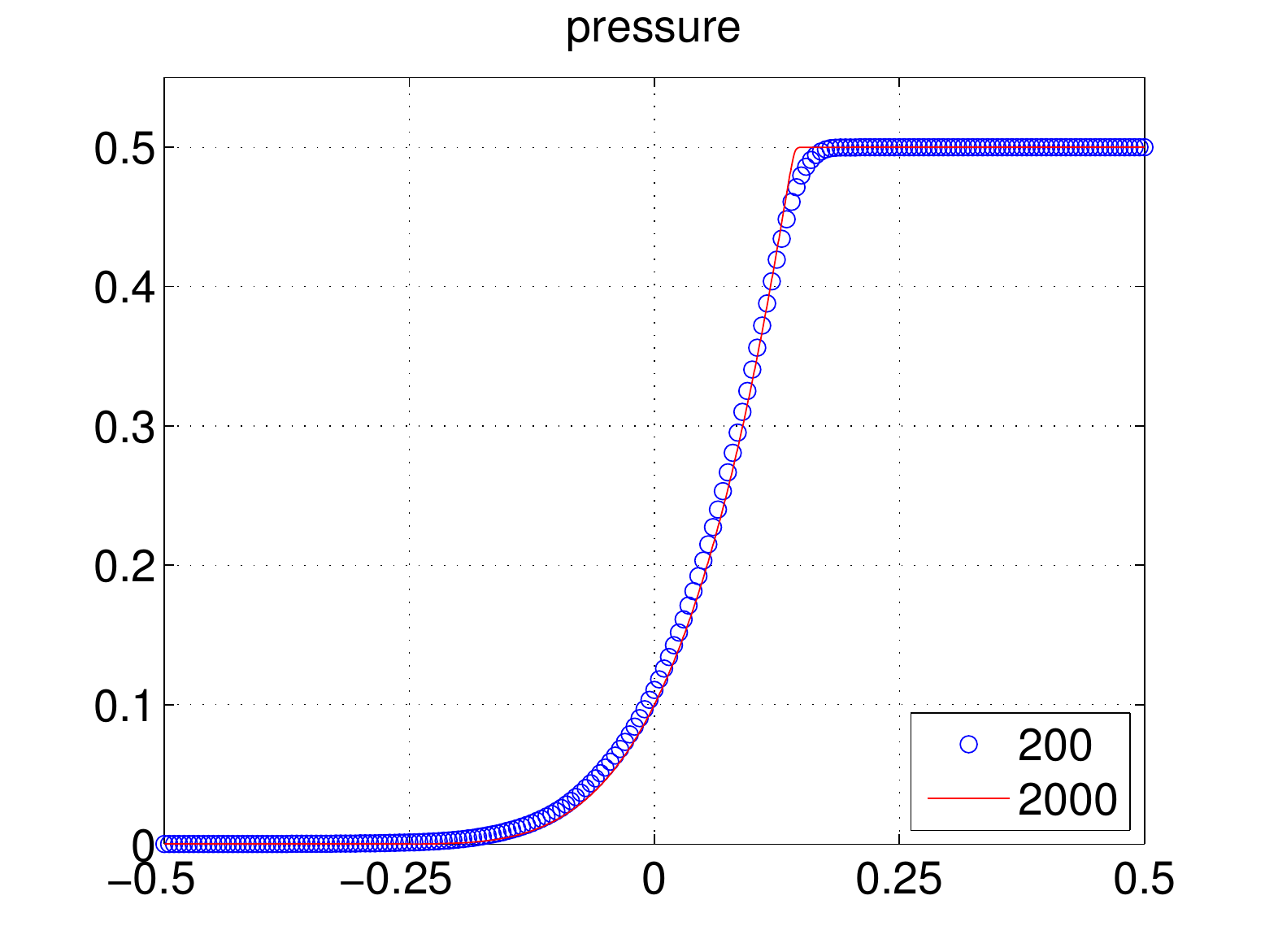} \\
\end{tabular}
      \caption{Vacuum shock tube problem. 
      Shown in these panels are plots at time $t=0.1$ of 
      (a) the density and
      (b) the thermal pressure. 
      The blue circle is a solution solved on a mesh with $N = 200$.
      The solid line is a highly resolved solution with $N = 2000$.     
      \label{vacuum}}
\end{center}
\end{figure}

We first consider a 1D vacuum shock tube problem.
This example is used to demonstrate our positivity-preserving MHD solver can handle very low density and pressure. 
The initial condition is:
\begin{align}
(\rho, u_x,u_y,u_z,p,B_x,B_y,B_z) =
\begin{cases}(10^{-12},0, 0, 0, 10^{-12}, 0,0,0) & \mbox{if} \quad x<0, \\ 
( 1, 0, 0, 0, 0.5, 0,1,0) & \mbox{if} \quad x>0. \end{cases}
\end{align}
It is similar to the vacuum shock tube problem in \cite{article:Wa09}.
The computational domain is $ [-0.5,0.5]$ and zero-order extrapolation boundary conditions are used.
Shown in Figure \ref{vacuum} are the density and pressure of the solution on a mesh with $N=200$ and the highly resolved solution with $N=2000$. 
We can observe the solution of low resolution and high resolution are in good agreements. 

\subsubsection{Torsional Alfv\'en wave pulse}
\begin{figure}
\begin{center}
\begin{tabular}{cc}
  	(a)\includegraphics[width=0.42\textwidth]{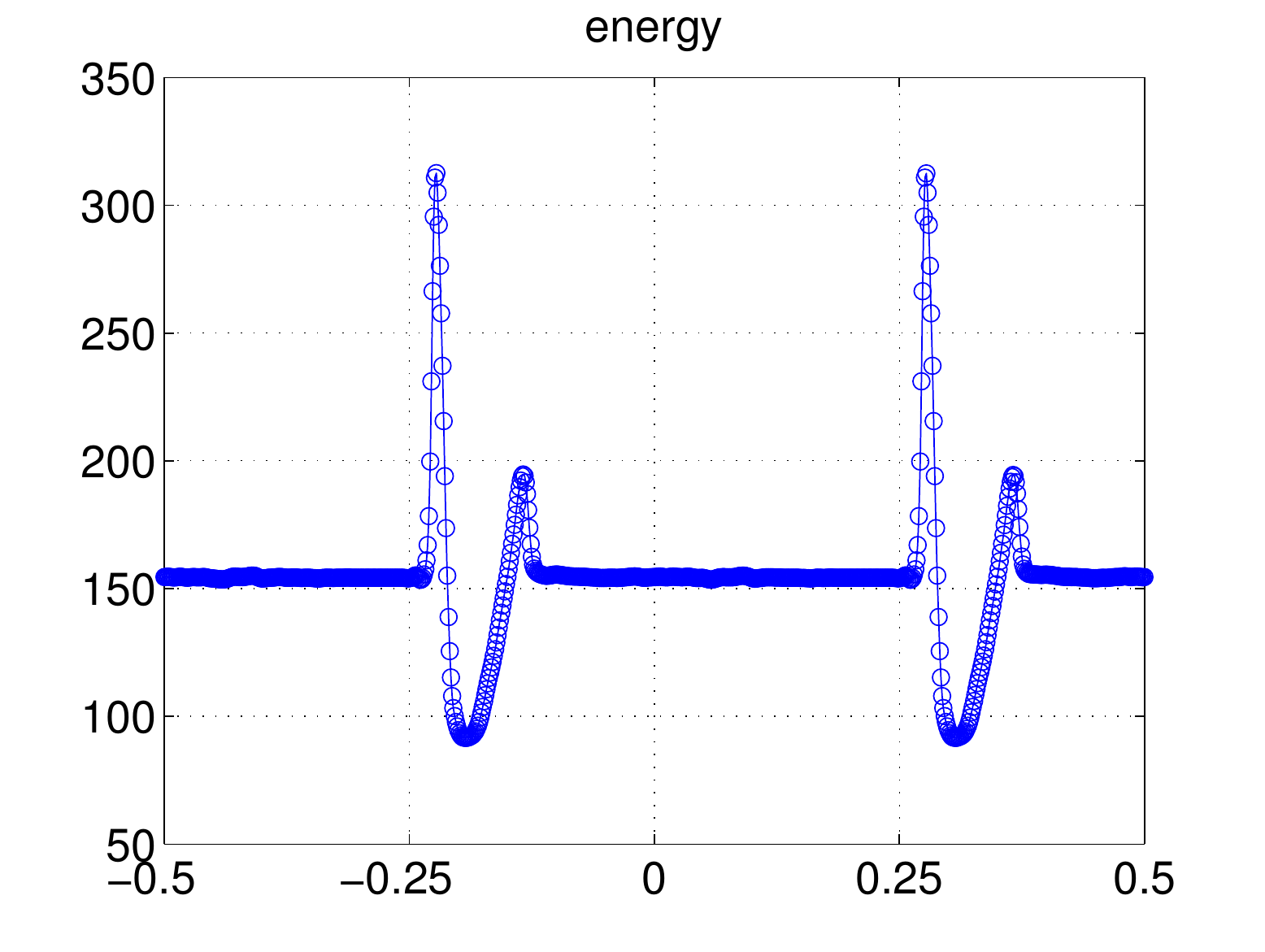} & 
	(b)\includegraphics[width=0.42\textwidth]{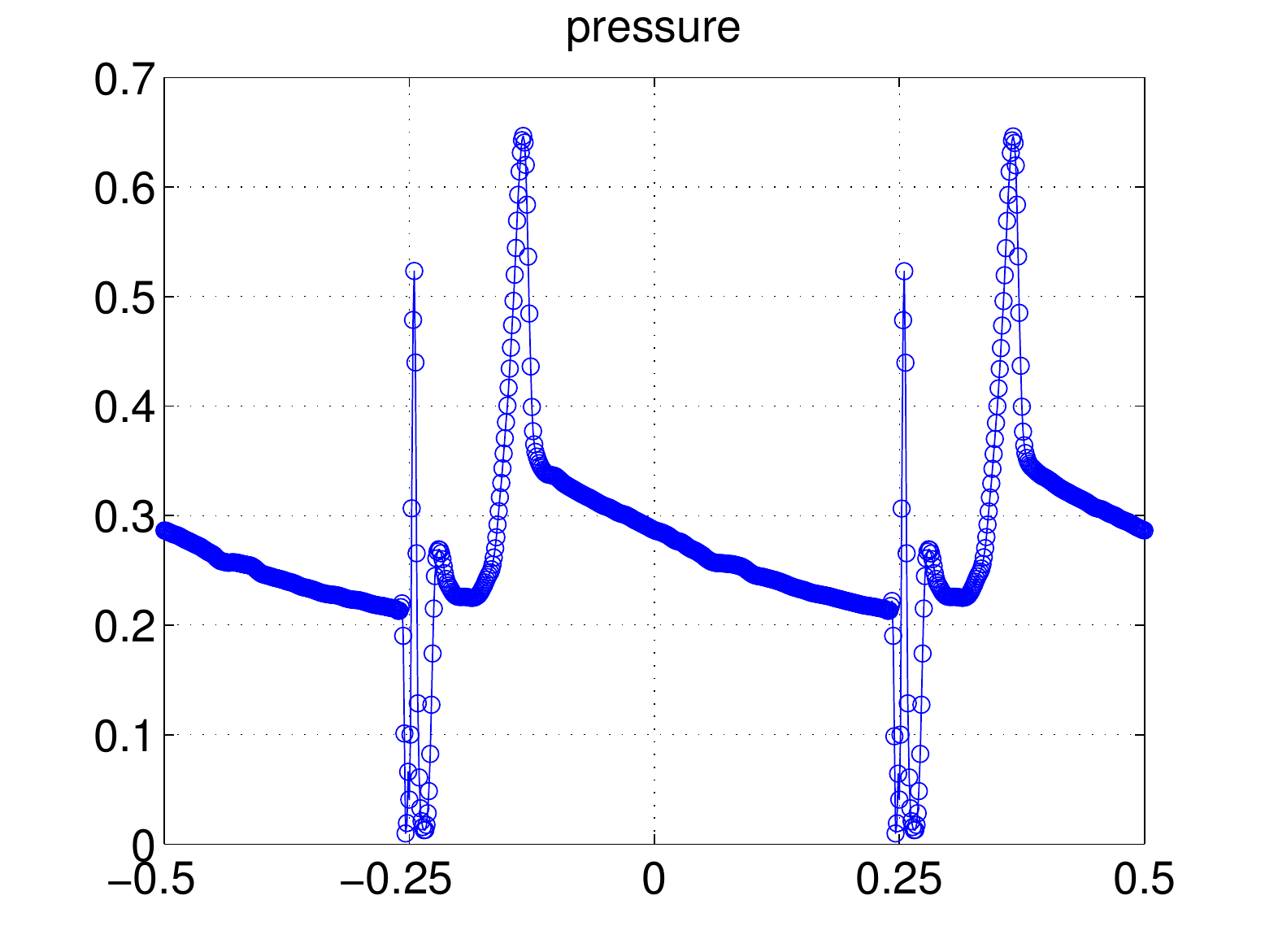}
\end{tabular}
      \caption{The torsional Alfv\'en wave pulse. 
      Shown in these panels are plots at time $t=0.156$ of 
      (a) the energy and
      (b) the thermal pressure.
      The solution was obtained on a mesh with $N=800$.     
      \label{alfven1}}
\end{center}
\end{figure}

\begin{figure}
\begin{center}
\begin{tabular}{cc}
  	(a)\includegraphics[width=0.42\textwidth]{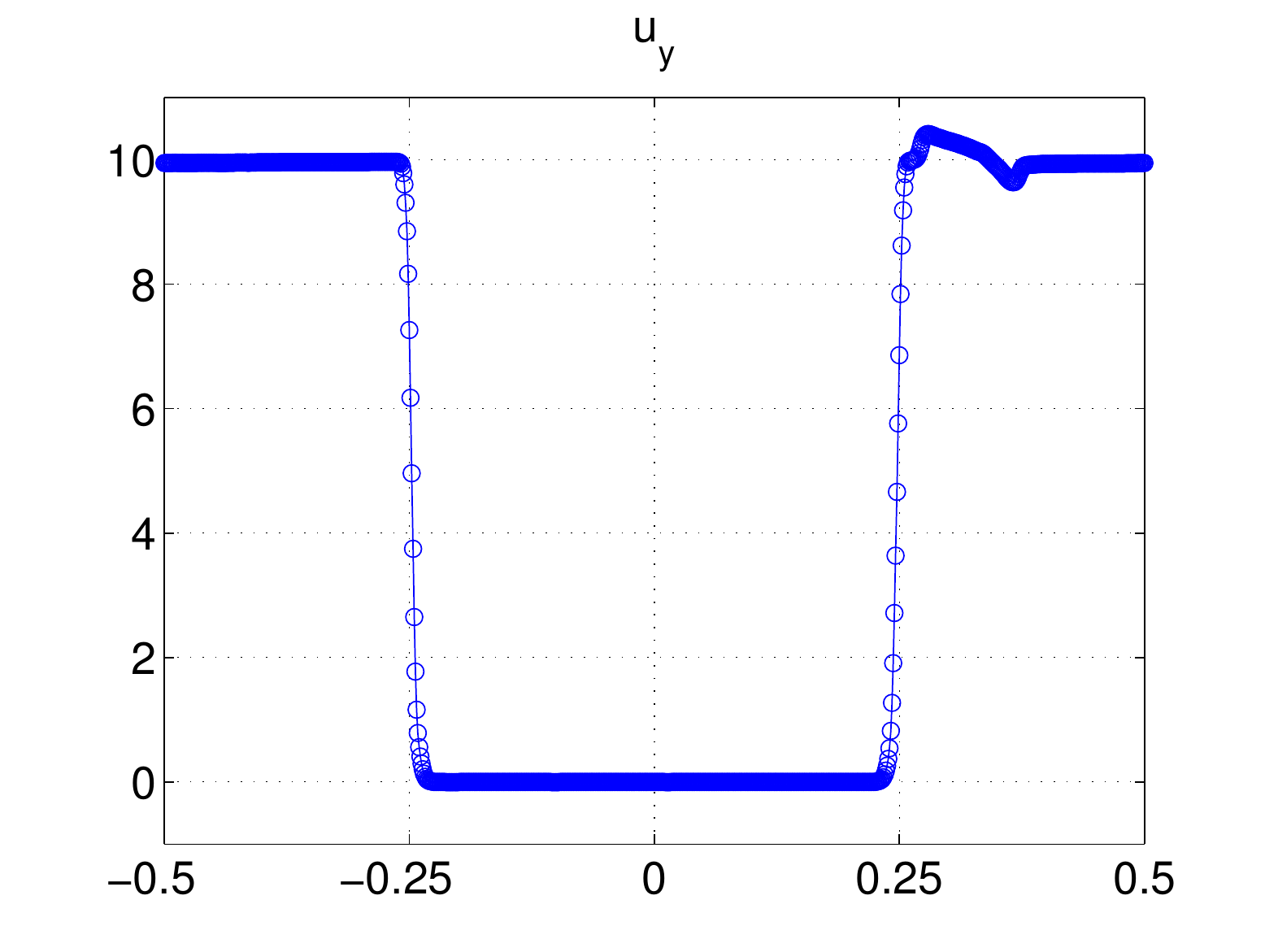} & 
	(b)\includegraphics[width=0.42\textwidth]{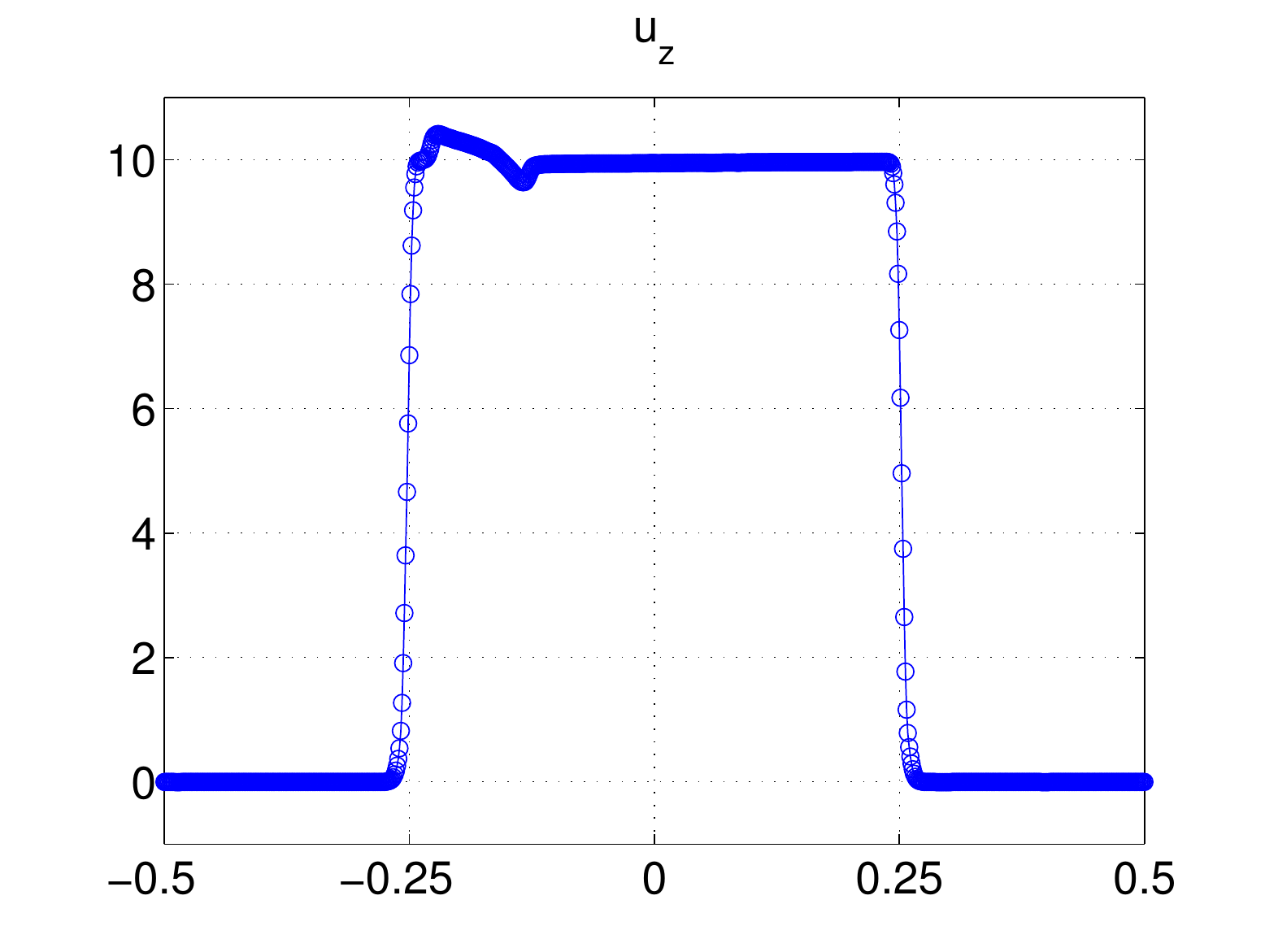} \\
	(c)\includegraphics[width=0.42\textwidth]{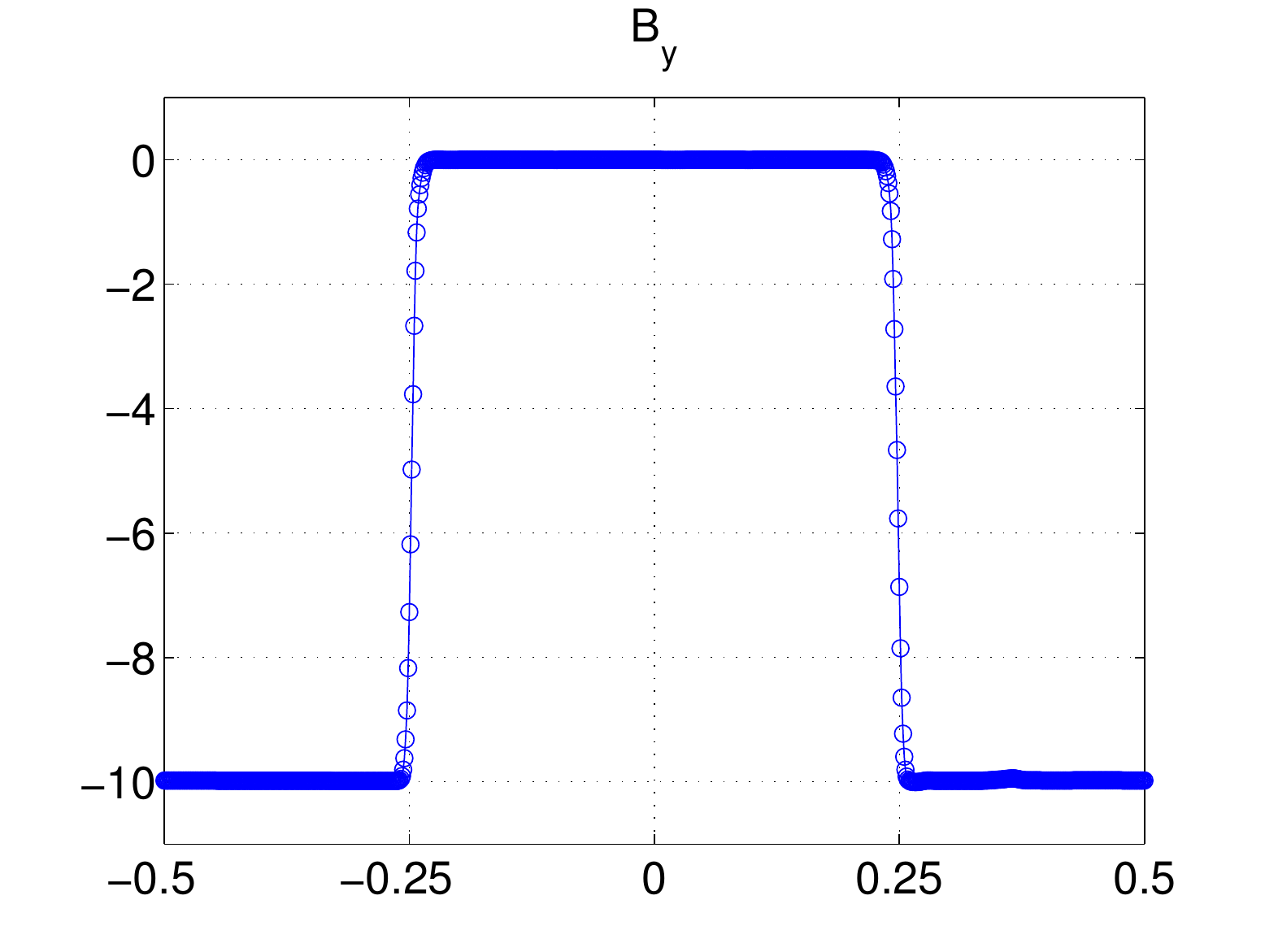} & 
	(d)\includegraphics[width=0.42\textwidth]{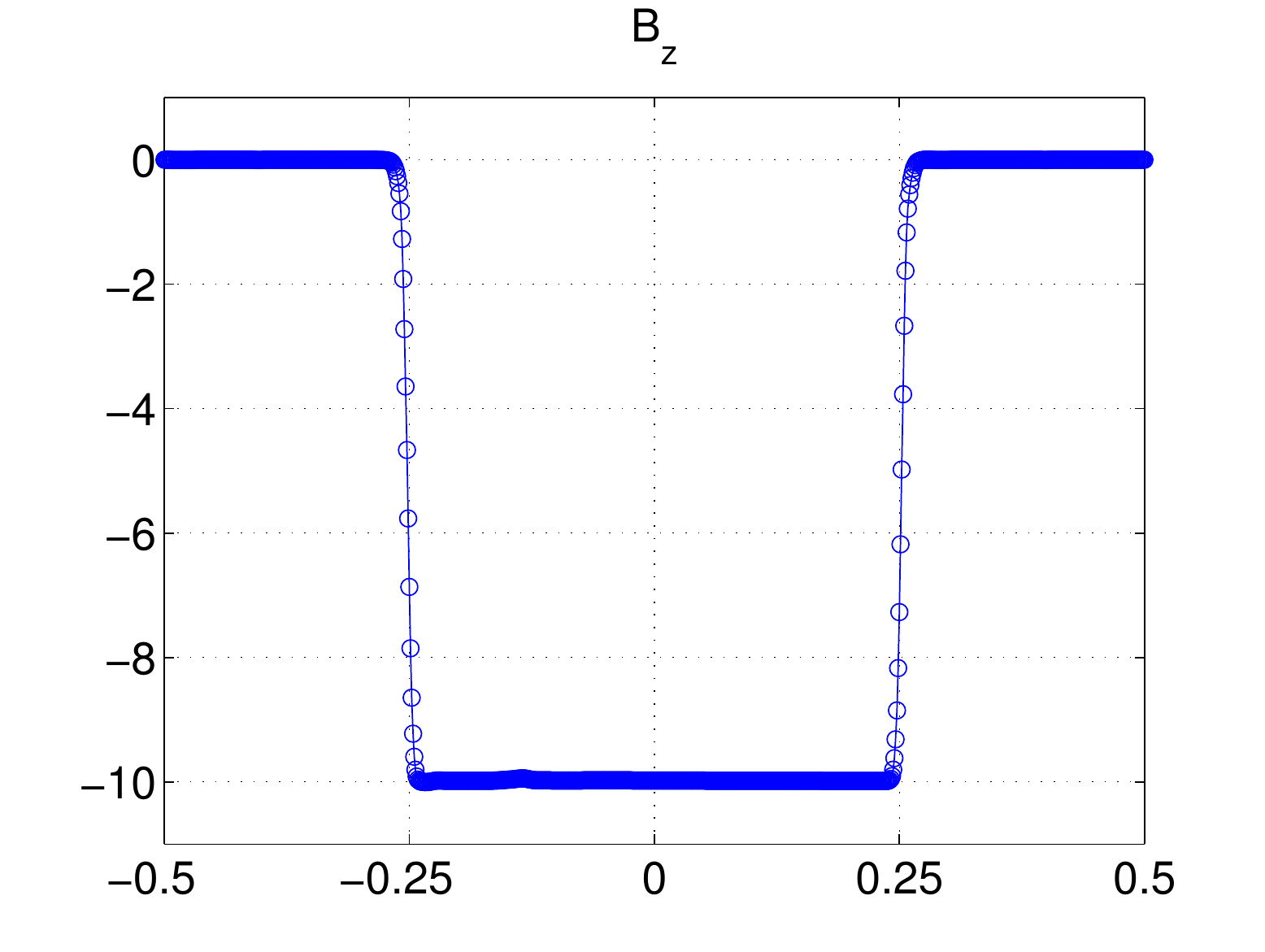}
\end{tabular}
      \caption{The torsional Alfv\'en wave pulse. 
      Shown in these panels are plots at time $t=0.156$ of 
      (a) $u_y$,
      (b) $u_z$,
      (c) $B_y$ and
      (d) $B_z$. 
      The solution was obtained on a mesh with $N=800$.    
      \label{alfven2}}
\end{center}
\end{figure}

We also consider the torsional Alfv\'en wave pulse problem \cite{article:BaSp99b, article:Cheng13}.
The initial condition is 
\begin{align}
(\rho, u_x,u_y,u_z,p,B_x,B_y,B_z) =
(1,10, 10 \cos \phi, 10 \sin \phi, 0.01, -10 \cos \phi, -10 \sin \phi,0),
\end{align}
where $\phi = \frac{\pi}{8} (\tanh(\frac{0.25+x}{\delta}+1))(\tanh(\frac{0.25-x}{\delta}+1))$
and $\delta = 0.005$. The computational domain is $[-0.5, 0.5]$ and periodic boundary conditions are used.
In this test problem, the initial pressure is so small 
that the problem is very sensitive to the dissipation introduced by numerical schemes. Further, the existence of a strong torsional Alfv\'en wave discontinuity makes the problem difficult to simulate. In the simulation without the  limiter, the base WENO-HCL introduced a negative pressure in a few time steps and the solutions become unphysical immediately. 
With the limiter, our scheme can simulate the problem stably and the numerical results at $t = 0.156$ are shown in Figure \ref{alfven1} and \ref{alfven2} with $N = 800$. 
Shown in the figures are plots of the energy, the thermal pressure,  $u_y$, $u_z$, $B_y$ and $B_z$. 
It is observed that our method successfully captures the two discontinuities and the results are comparable with those in \cite{article:BaSp99b, article:LiXuYa11}.
However, small bumps can still be observed around one of the discontinuities of both $u_y$ and $u_z$.
The authors in \cite{article:BaSp99b} pointed out this is because the MHD solver introduced too much numerical dissipation to keep the pressure positive.
The primary reason is the Riemann solver around the discontinuities is not selective enough.

\subsection{Test cases in multi-D}
In this subsection, we consider several 2D and 3D examples to demonstrate the accuracy and efficiency of our positivity-preserving multi-D MHD solver in CT framework.
In the following tests, we implement fourth-order WENO-CT2D and WENO-CT3D schemes 
as the MHD solver, to which we apply our positivity-preserving limiter.
Unless otherwise stated, we use {\bf Option 2} for the multi-D simulation. 


\subsubsection{Smooth vortex test in MHD}
We consider the smooth vortex problem with non-zero magnetic field
to demonstrate the scheme can maintain the designed accuracy within the CT framework.
We consider a modification of the smooth vortex problem considered in \cite{article:Ba04, article:LiXu12, xu2011divergence}.
The initial condition is a mean flow
\begin{align}
(\rho, u_x, u_y, u_z, p, B_x, B_y, B_z) = (1, 1, 1, 0, 1, 0, 0, 0),
\end{align}
with perturbations on $u_x$, $u_y$, $B_x$, $B_y$ and $p$:
\begin{gather*}
(\delta u_x, \delta u_y) = \frac{\kappa}{2\pi}e^{0.5(1-r^2)}(-{y},{x}), \qquad
(\delta B_x, \delta B_y) = \frac{\mu}{2\pi}e^{0.5(1-r^2)}(-{y},{x}), \\
\delta p = \frac{\mu^2(1-r^2) - \kappa^2}{8\pi^2}e^{1-r^2}.
\end{gather*}
The magnetic potential is initialized as
\begin{gather*}
A_z = \frac{\mu}{2\pi} e^{0.5(1-r^2)}
\end{gather*}
Here $r^2={x}^2 + {y}^2$. 

We set the vortex strength $\mu=5.389489439$ and $\kappa = \sqrt{2} \mu$ such that the lowest 
pressure in the center of the vortex is $5.3 \times 10^{-12}$.
Similar to \cite{article:LiXu12}, we use computational domain $(x,y) \in [-10,10]\times[-10,10]$ 
such that the error from the boundary conditions will not influence the overall convergence study.
The periodic boundary condition are used on all sides.
Because fourth-order CT steps are used, the overall scheme is fourth-order accuracy.

The $L_1$-errors and $L_\infty$-errors of the velocity and magnetic field for $t = 0.05$ are shown in Tables~\ref{tab:2dmhd}, in which one can conclude the proposed positivity-preserving scheme can maintain fourth-order accuracy as expected. 
We remark that negative pressure is observed on meshes coarser than $320 \times 320$
when the limiter is not applied.


\begin{table}
\begin{center}
\begin{Large}
    \caption{Accuracy test of the 2D vortex evolution in MHD. Shown are the $L_{1}$-errors and $L_{\infty}$-errors at time $t = 0.05$ of the density as computed by the positivity-preserving WENO-CT2D scheme at various grid resolutions. 
    The solutions converge at fourth-order accuracy.   \label{tab:2dmhd}}
    \end{Large}

\begin{tabular}{|c || c c c c|| c c c c| }
  \hline
{\normalsize{ Mesh}} & \multicolumn{4}{c||}{$u_x$} & \multicolumn{4}{c|}{$u_y$} \\ \cline{2-5} \cline{6-9}
 & {\normalsize  $L_{1}$-Error}  &  {\normalsize  Order} & {\normalsize $L_{\infty}$-Error } & {\normalsize Order }
 & {\normalsize  $L_{1}$-Error}  &  {\normalsize  Order} & {\normalsize $L_{\infty}$-Error } & {\normalsize Order }
 \\ \hline \hline
{\normalsize $40 \times 40$}    & 7.38E-04	&	-	&	1.79E-02	&	-	
								& 8.03E-04	&	-	&	1.94E-02	&	-\\
								
{\normalsize $80 \times 80$}  	&7.20E-05	&	3.35	&	4.33E-03& 	2.05	
								&7.36E-05	&	3.45	&	5.22E-03& 	1.90\\
								
{\normalsize $160 \times 160$}  &3.46E-06	&	4.38	&	1.92E-04&	4.49	
								 &3.72E-06	&	4.31	&	2.18E-04&	4.58\\
								
{\normalsize $320 \times 320$} 	&1.80E-07	&	4.27	&	1.42E-05&	3.76	
								&1.96E-07	&	4.25	&	1.64E-05&	3.73\\  \hline
 \end{tabular}
 
 \vspace{.2in}
 
\begin{tabular}{|c || c c c c|| c c c c |}
  \hline
{\normalsize{ Mesh}} & \multicolumn{4}{c||}{$B_x$} & \multicolumn{4}{c|}{$B_y$} \\ \cline{2-5} \cline{6-9}
 & {\normalsize  $L_{1}$-Error}  &  {\normalsize  Order} & {\normalsize $L_{\infty}$-Error } & {\normalsize Order }
 & {\normalsize  $L_{1}$-Error}  &  {\normalsize  Order} & {\normalsize $L_{\infty}$-Error } & {\normalsize Order }
 \\ \hline \hline
{\normalsize $40 \times 40$}    & 1.02E-03	&	-	&	1.49E-02	&	-	
								& 1.04E-03	&	-	&	1.57E-02	&	-\\
								
{\normalsize $80 \times 80$}  	&7.73E-05	&	3.73	&	1.27E-03& 	3.56	
								&7.73E-05	&	3.75	&	1.16E-03& 	3.77\\
								
{\normalsize $160 \times 160$}  &4.75E-06	&	4.03	&	8.25E-05&	3.94	
								 &4.74E-06	&	4.03	&	7.16E-05&	4.01\\
								
{\normalsize $320 \times 320$} 	&2.85E-07	&	4.06	&	7.66E-06&	3.43	
								&2.84E-07	&	4.06	&	6.36E-06&	3.49\\  \hline
 \end{tabular}
\end{center}
\end{table}

\subsubsection{Rotated vacuum shock tube problem}
\label{section:2drs}
In this example, we consider the vacuum shock tube problem rotated by an angle of $\alpha$ in a 2D domain. The initial conditions in this case are
\begin{align}
(\rho, u_{\perp},u_{\parallel},u_z,p,B_{\perp},B_{\parallel},B_z) =
\begin{cases}(10^{-12},0, 0, 0, 10^{-12}, 0,0,0) & \mbox{if} \quad \xi<0, \\ 
( 1, 0, 0, 0, 0.5, 0,1,0) & \mbox{if} \quad \xi>0. \end{cases}
\end{align}
where 
\begin{equation}
\xi = x \cos \alpha + y \sin \alpha \quad \text{and} \quad \eta = -x \sin \alpha + y \cos \alpha,
\end{equation}
where $u_{\perp}$ and $B_{\perp}$ are perpendicular to the shock interface, and $u_{\parallel}$ and $B_{\parallel}$  are parallel to the shock interface. The magnetic potential is initialized as
\begin{align}
A_z(0,\xi) =
\begin{cases} 0 & \mbox{if} \quad \xi \le 0, \\ 
- \xi & \mbox{if}  \quad \xi \ge 0. \end{cases}
\end{align}
We solve this problem by the positivity-preserving WENO-CT2D scheme on the computational domain $(x,y) \in [-0.6,0.6] \times [-0.25,0.25]$  with a $240 \times 100$ mesh. $\alpha = \tan^{-1}(0.5)$.
Zero-order extrapolation boundary conditions are used on the left and right boundaries.
On the top and bottom boundaries, all the quantities are set to describe the exact motion of the shock.

\begin{figure}
\begin{center}
\begin{tabular}{cc}
  	(a)\includegraphics[width=0.42\textwidth]{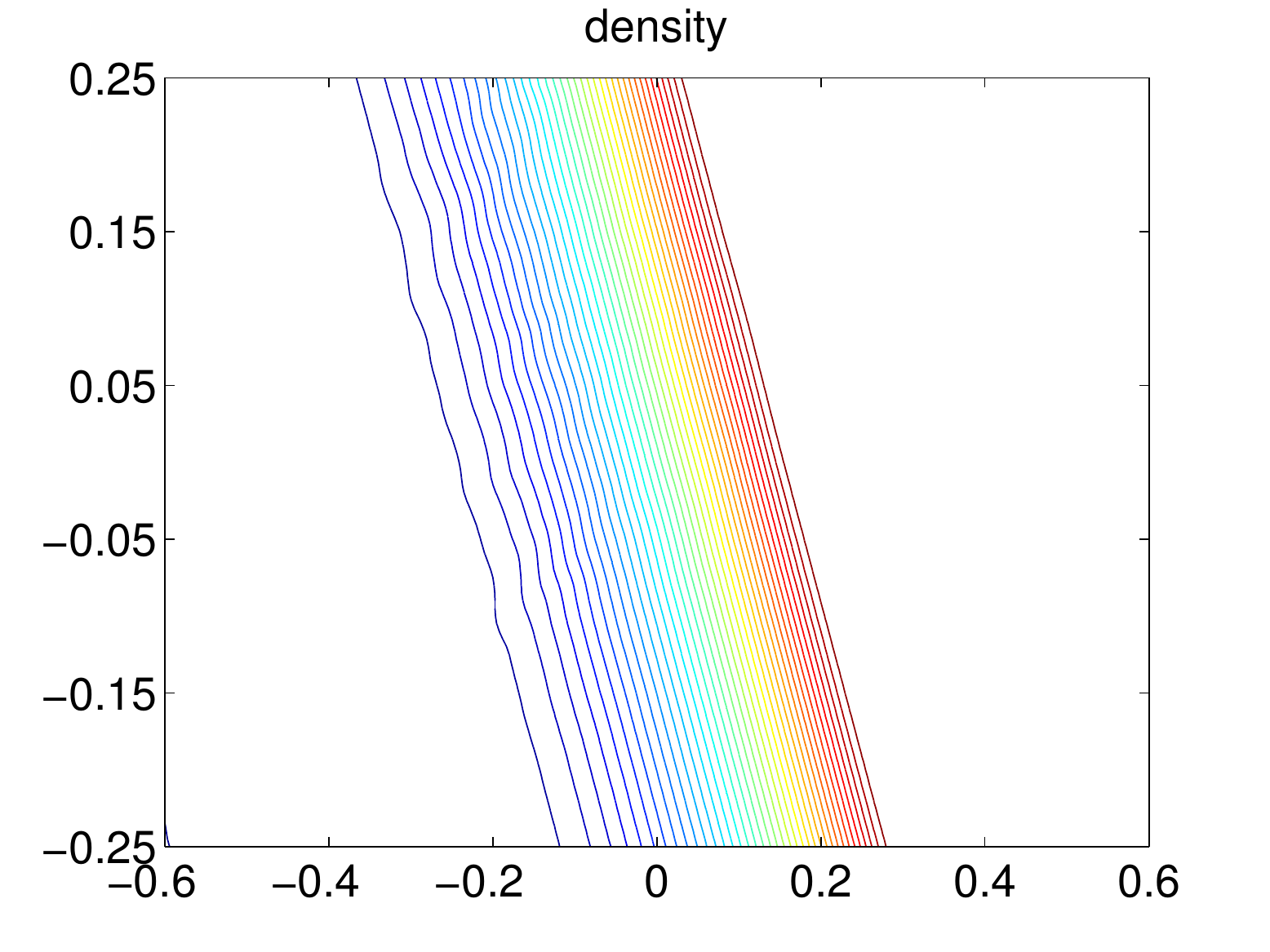} & 
	(b)\includegraphics[width=0.42\textwidth]{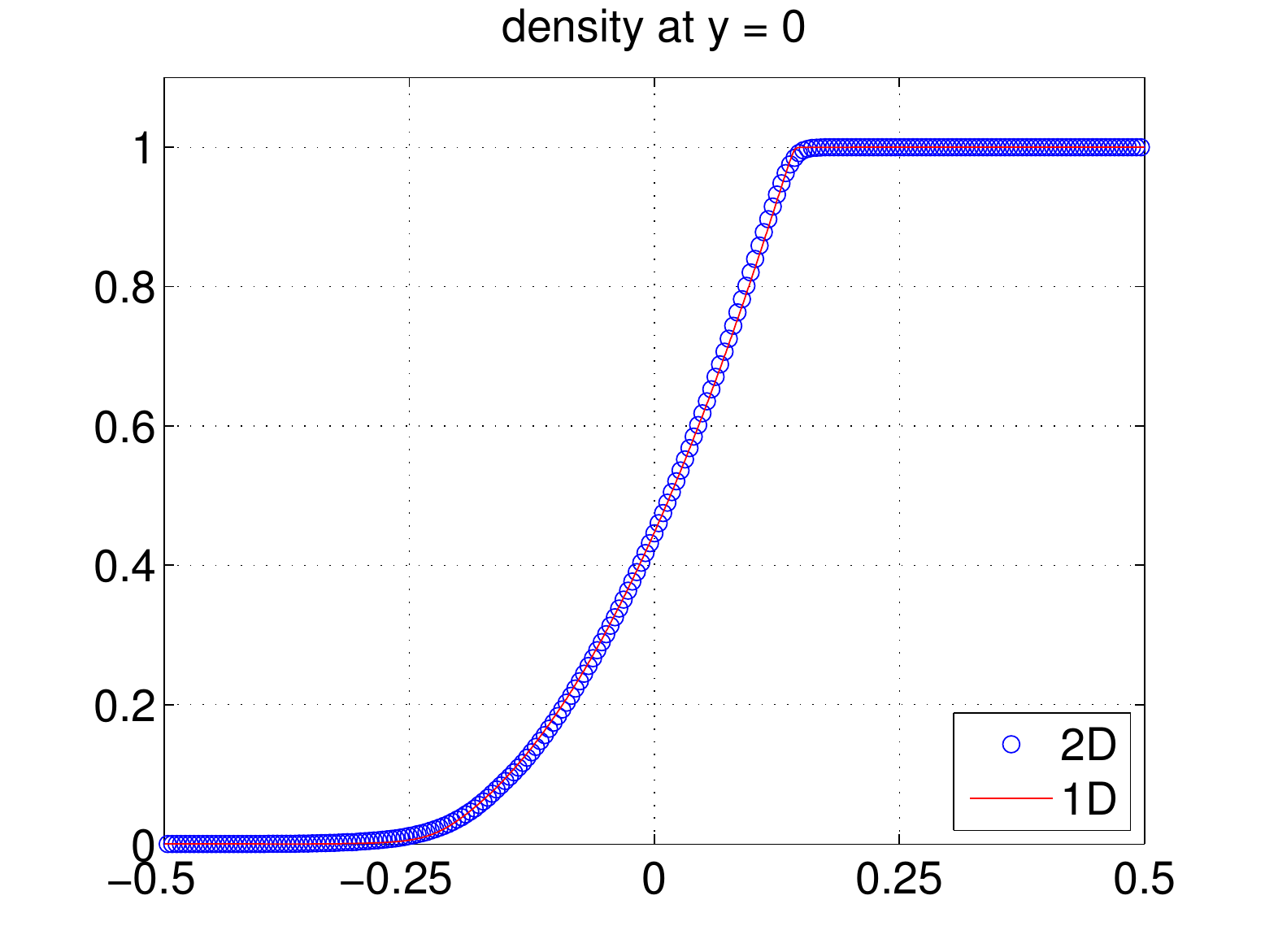} \\
	(c)\includegraphics[width=0.42\textwidth]{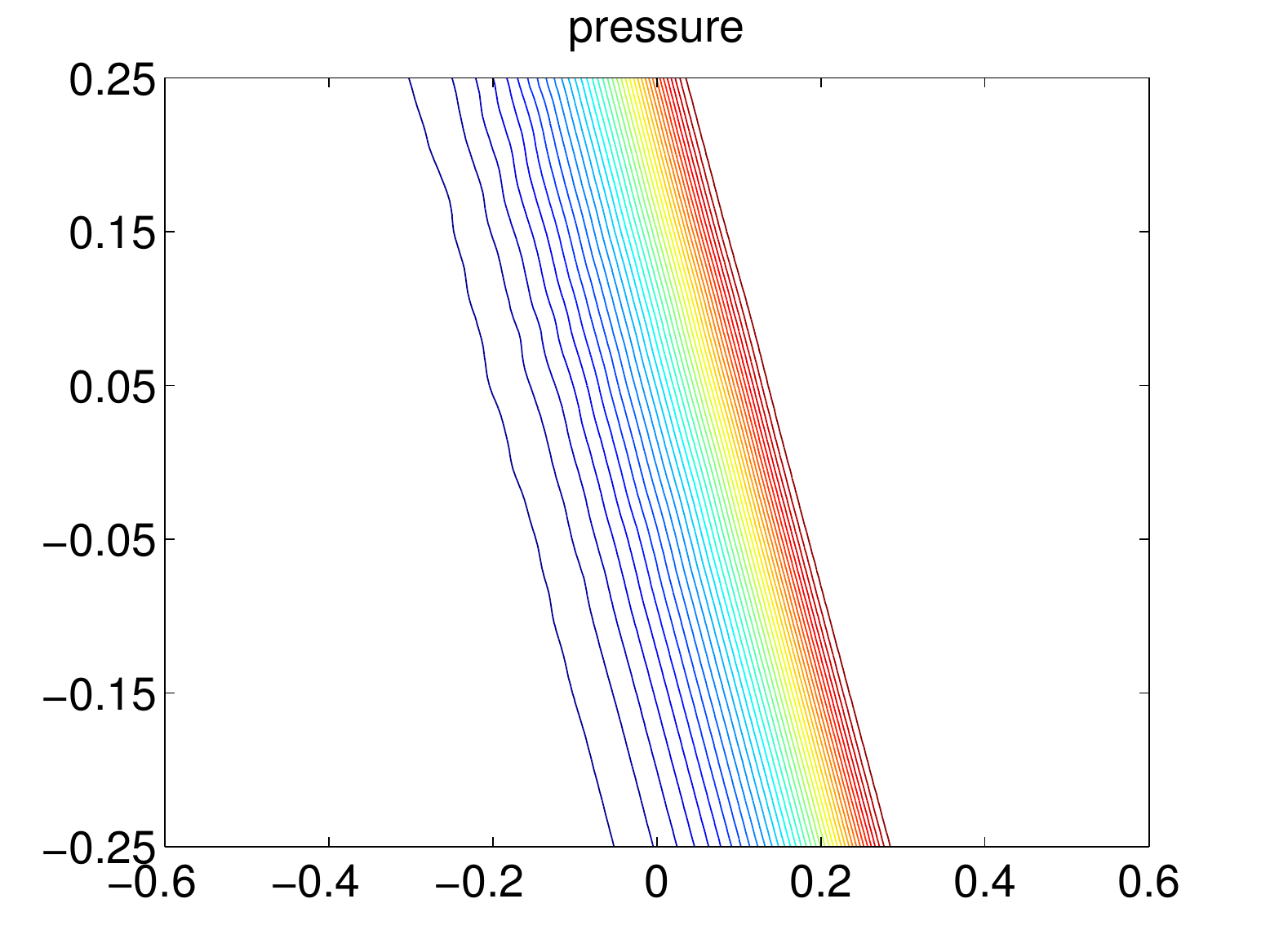} & 
	(d)\includegraphics[width=0.42\textwidth]{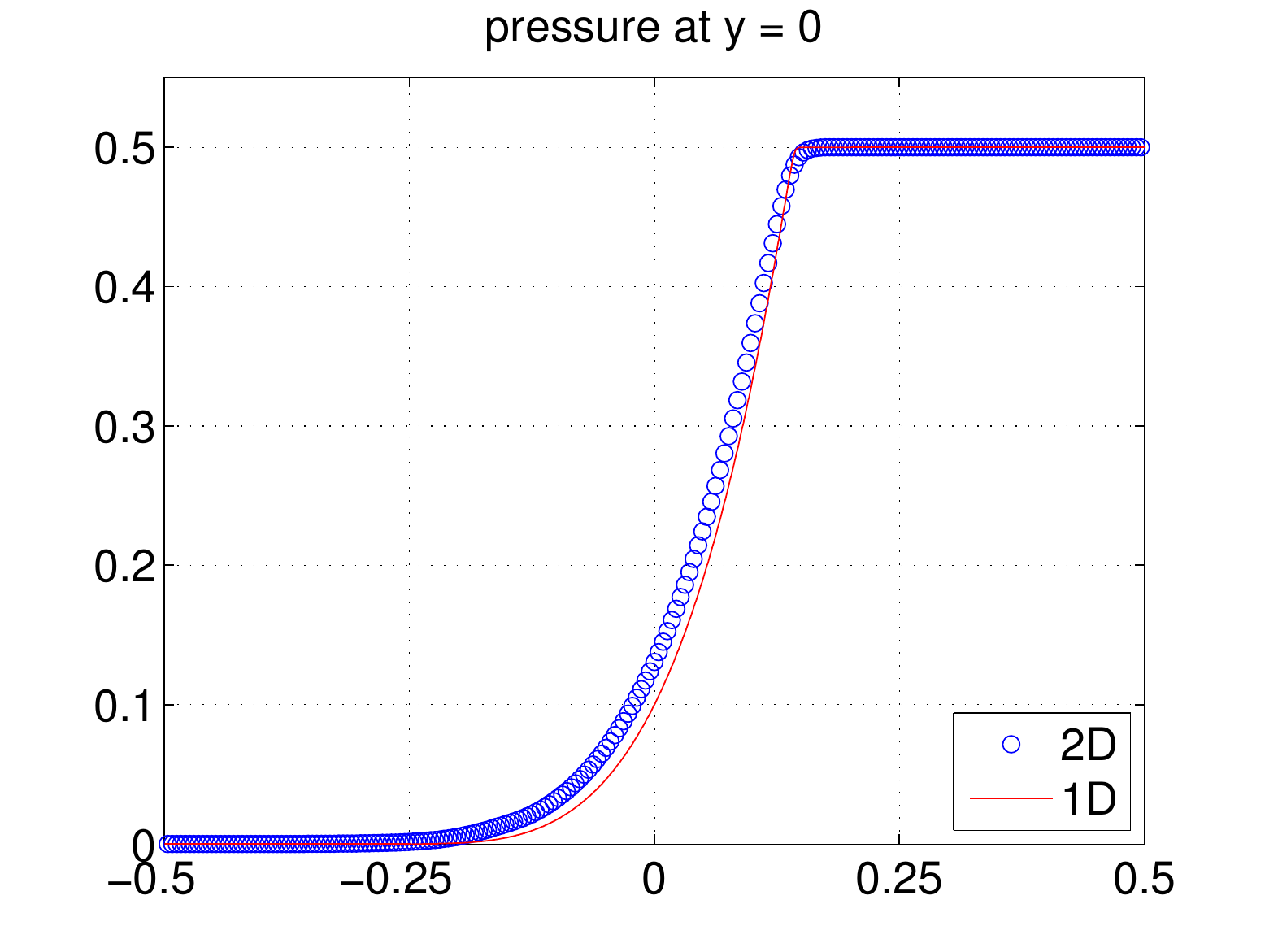}
\end{tabular}
      \caption{Rotated vacuum shock tube problem. 
      Shown in these panels are plots at time $t=0.1$ of 
      (a) the density,
      (b) the density cut at y = 0,
      (c) the pressure and
      (d) the pressure cut at y = 0. 
      40 equally spaced contours are used for the contour plots. The solid lines in (b) and (d)
      are 1D highly resolved solutions.
      The solution was obtained on a $240 \times 100$ mesh.     
      \label{vacuum2d}}
\end{center}
\end{figure}

The solutions are plotted in Figure \ref{vacuum2d}, where 1D cut of density and pressure at $y = 0$ is also plotted
to compare with the 1D highly resolved results. 
We clearly observe that the 2D solution is consistent with the 1D solution.
Without the limiter, negative density and pressure are observed in numerical solutions, 
which quickly leads to blow-up of the numerical simulation.

\subsubsection{2D blast problem}
In the blast wave problem, a strong fast magnetosonic shock formulates and propagates into the low-$\beta$  plasma background, which will likely lead to negative density or pressure in numerical solutions.
In this subsection, we first investigate a 2D version of the problem \cite{article:Bal09-ader, article:BaSp99a, article:LiXuYa11}.
The computational domain is $(x,y) \in [-0.5,0.5] \times [-0.5, 0.5]$ with outflow boundary conditions on all the four sides. 

\begin{figure}
\begin{center}
\begin{tabular}{cc}
  	(a)\includegraphics[width=0.42\textwidth]{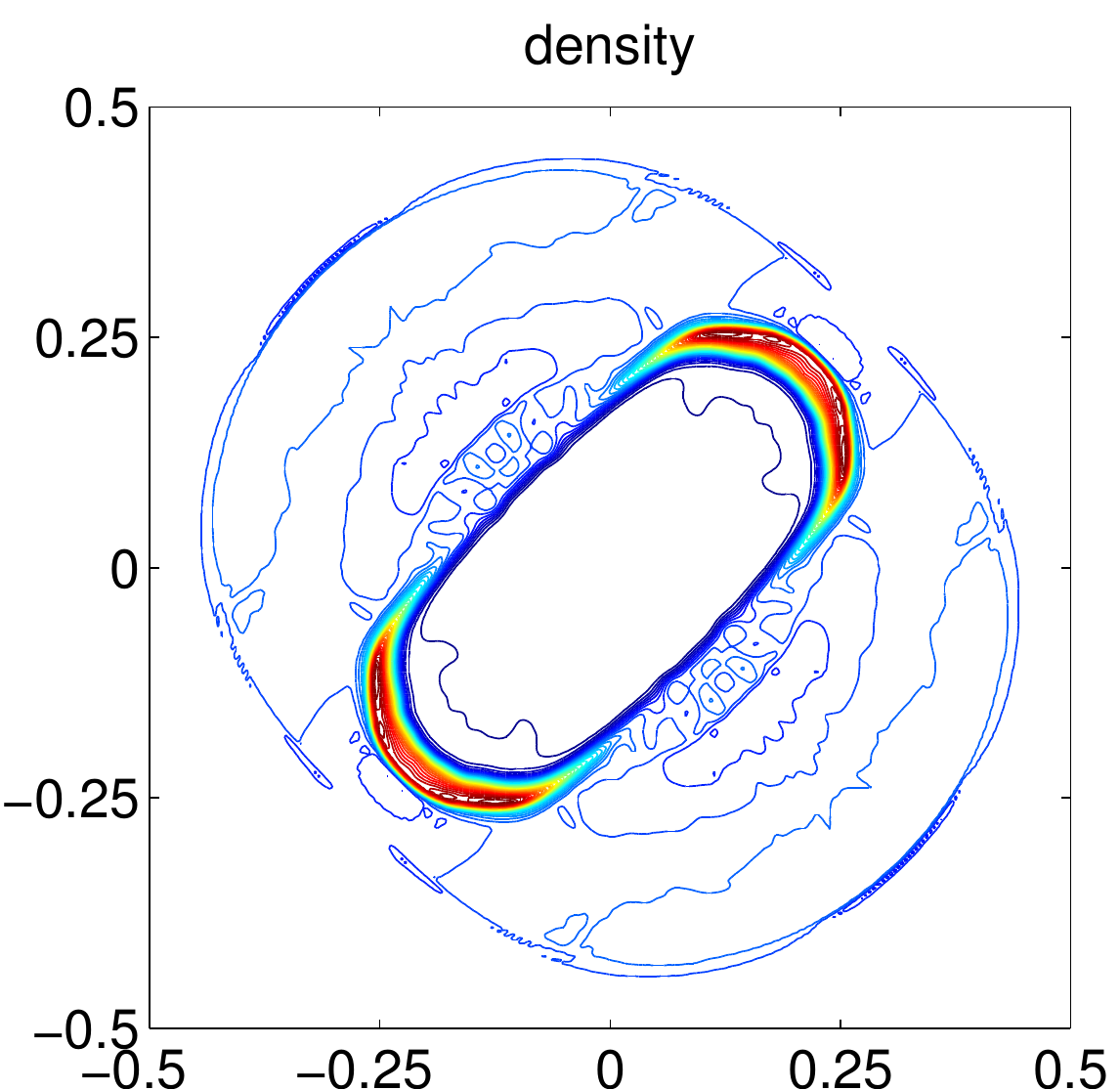} & 
	(b)\includegraphics[width=0.42\textwidth]{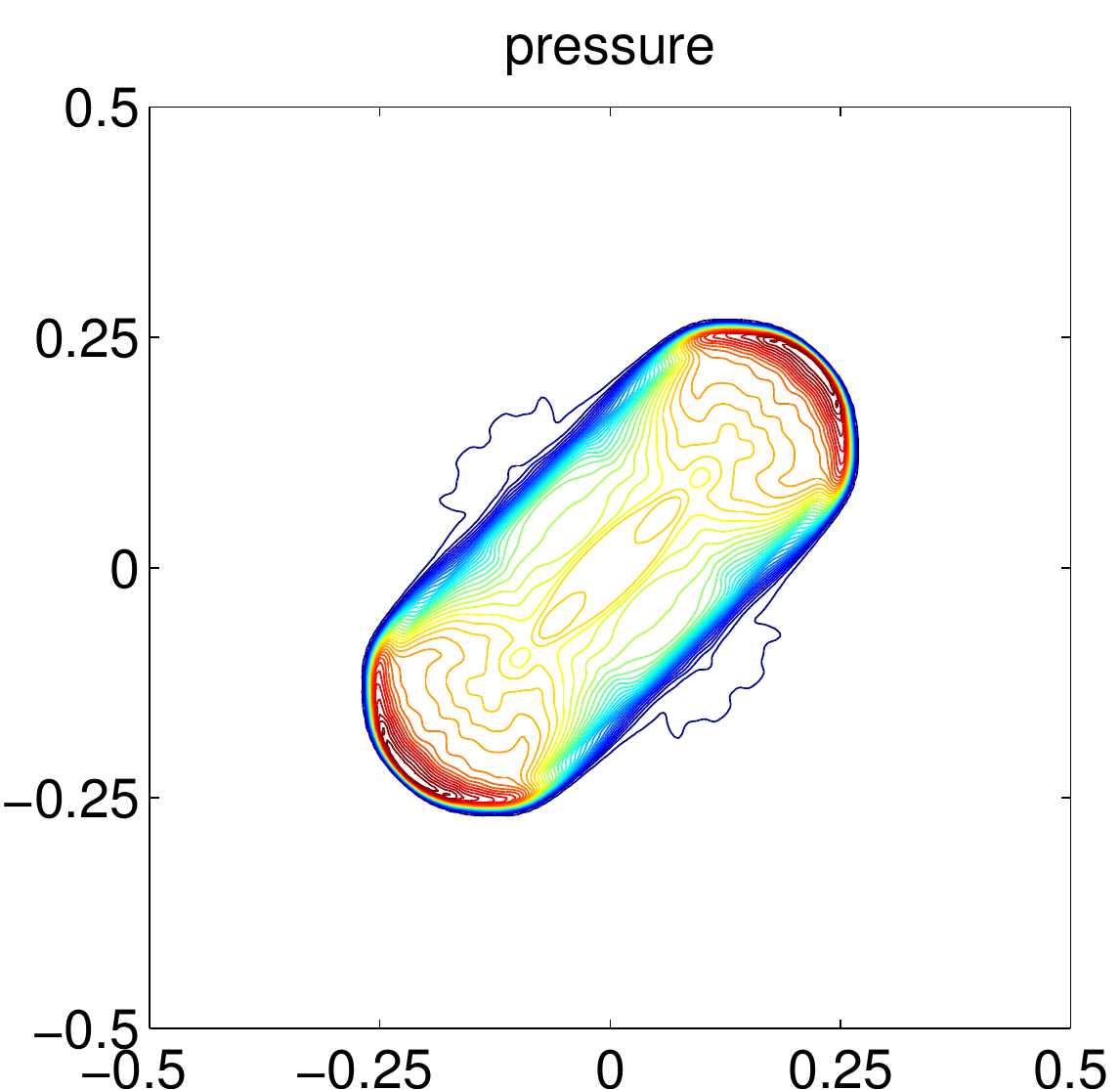} \\
	(c)\includegraphics[width=0.42\textwidth]{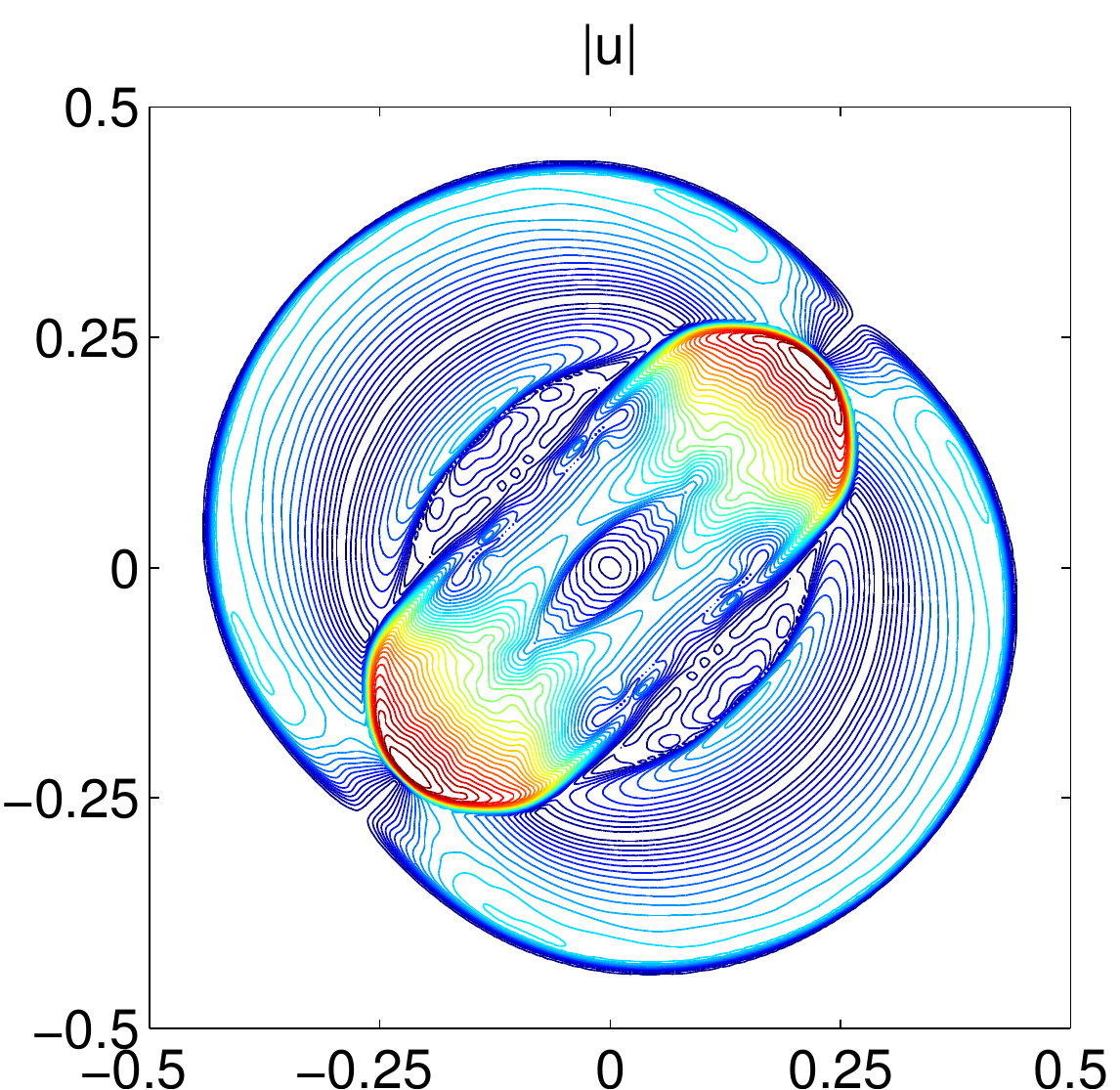} & 
	(d)\includegraphics[width=0.42\textwidth]{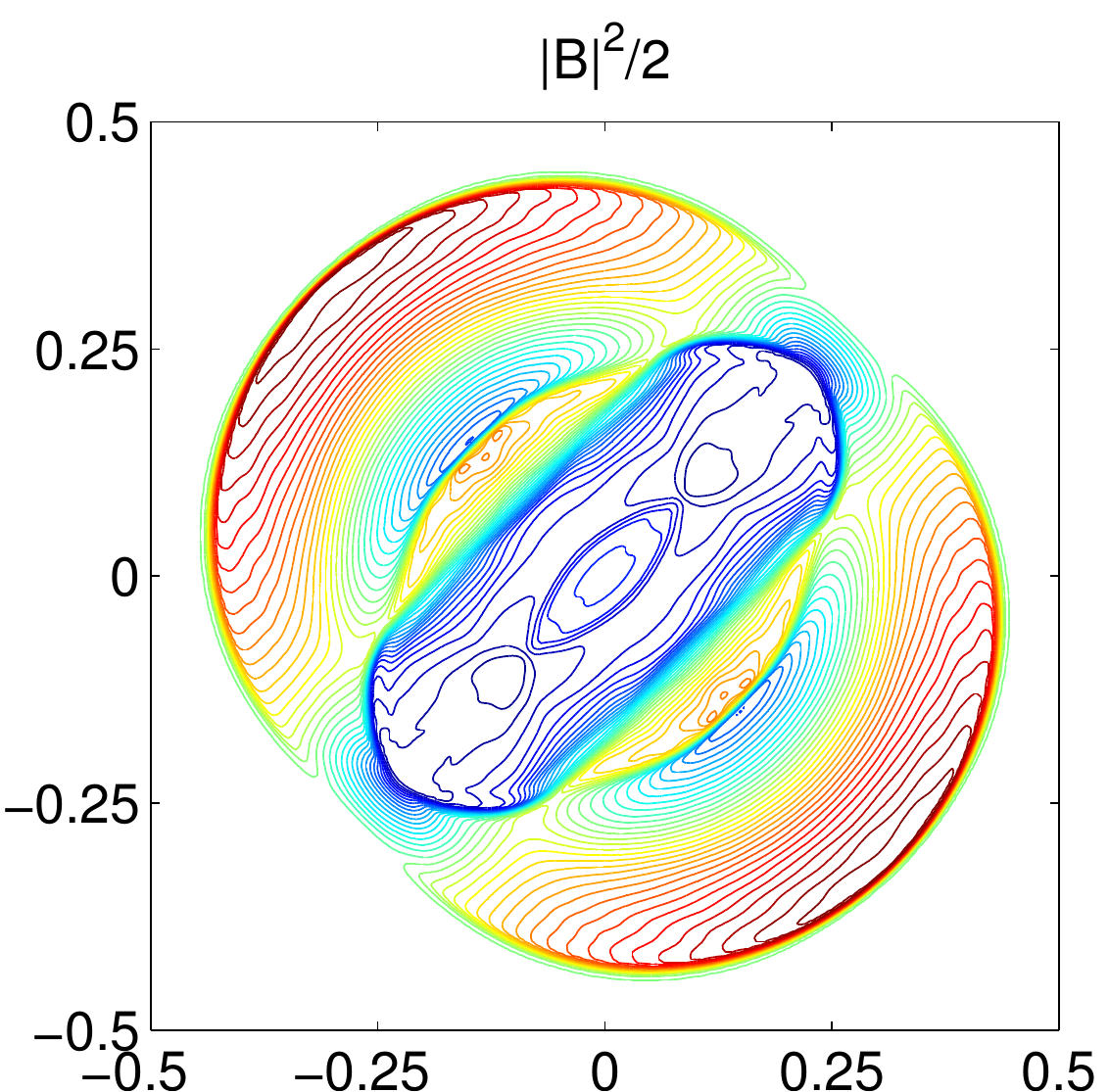}
\end{tabular}
      \caption{2D blast problem. 
      Shown in these panels are plots at time $t=0.01$ of 
      (a) the density,
      (b) the thermal pressure,
      (c) the norm of velocity and
      (d) the magnetic pressure. 
      40 equally spaced contours are used for each plot.
      The solution was obtained on a $256 \times 256$ mesh by positivity-preserving WENO-CT scheme.    
      \label{blast2}}
\end{center}
\end{figure}

The initial conditions of the problem consist of an initial background:
\begin{equation}
(\rho, u_x, u_y, u_z, p, B_x, B_y, B_z) = (1, 0, 0, 0, 0.1, 100/\sqrt{2 \pi},100/\sqrt{2 \pi},0),
\end{equation}
and a circular pressure pulse $p = 1000$ within a radius $r=0.1$ from the center of the domain. The initial scalar magnetic potential is simply given by
\begin{align}
A_z =  100/\sqrt{2 \pi} y - 100/\sqrt{2 \pi} x.
\end{align}
The solution is computed on a $256 \times 256$ mesh. Shown in Figure~\ref{blast2} are plots of the solutions. 
The solution shows good agreement with those in \cite{article:Bal09-ader,article:LiXuYa11}.

In Table~\ref{tab:comp}, we use this example to compare four different schemes, WENO-HCL, WENO-CT-OP1, WENO-CT-OP2 and PP-WENO-CT-OP2. 
Here WENO-HCL is referred to the base WENO-HCL scheme without CT or the limiter.
WENO-CT-OP1 is referred to the WENO-CT2D scheme choosing {\bf Option 1} without the limiter. 
WENO-CT-OP2 is referred to the WENO-CT2D scheme choosing {\bf Option 2} without the limiter. 
Finally, PP-WENO-CT-OP2 is referred to the positivity-preserving WENO-CT2D scheme with {\bf Option 2} chosen.

From Table~\ref{tab:comp}, we observe that the base WENO-HCL scheme is unstable for the resolution $150 \times 150$ and becomes stable in the higher resolutions.
WENO-CT-OP1 is unstable for each resolution and applying the positivity-preserving limiter will not be able to stabilize this because the negative pressure is from the correction step of the magnetic field.
WENO-CT-OP2 is stable in the lower resolution but become unstable for the resolution $256 \times 256$, and negative pressure is observed in all the resolutions.
Finally, the positivity-preserving WENO-CT scheme is stable for all the resolutions.
From those results, it is very clear that the positivity-preserving WENO-CT scheme is the most stable methods in those four methods.

\begin{figure}
\begin{center}
\begin{tabular}{cc}
  	(a)\includegraphics[width=0.42\textwidth]{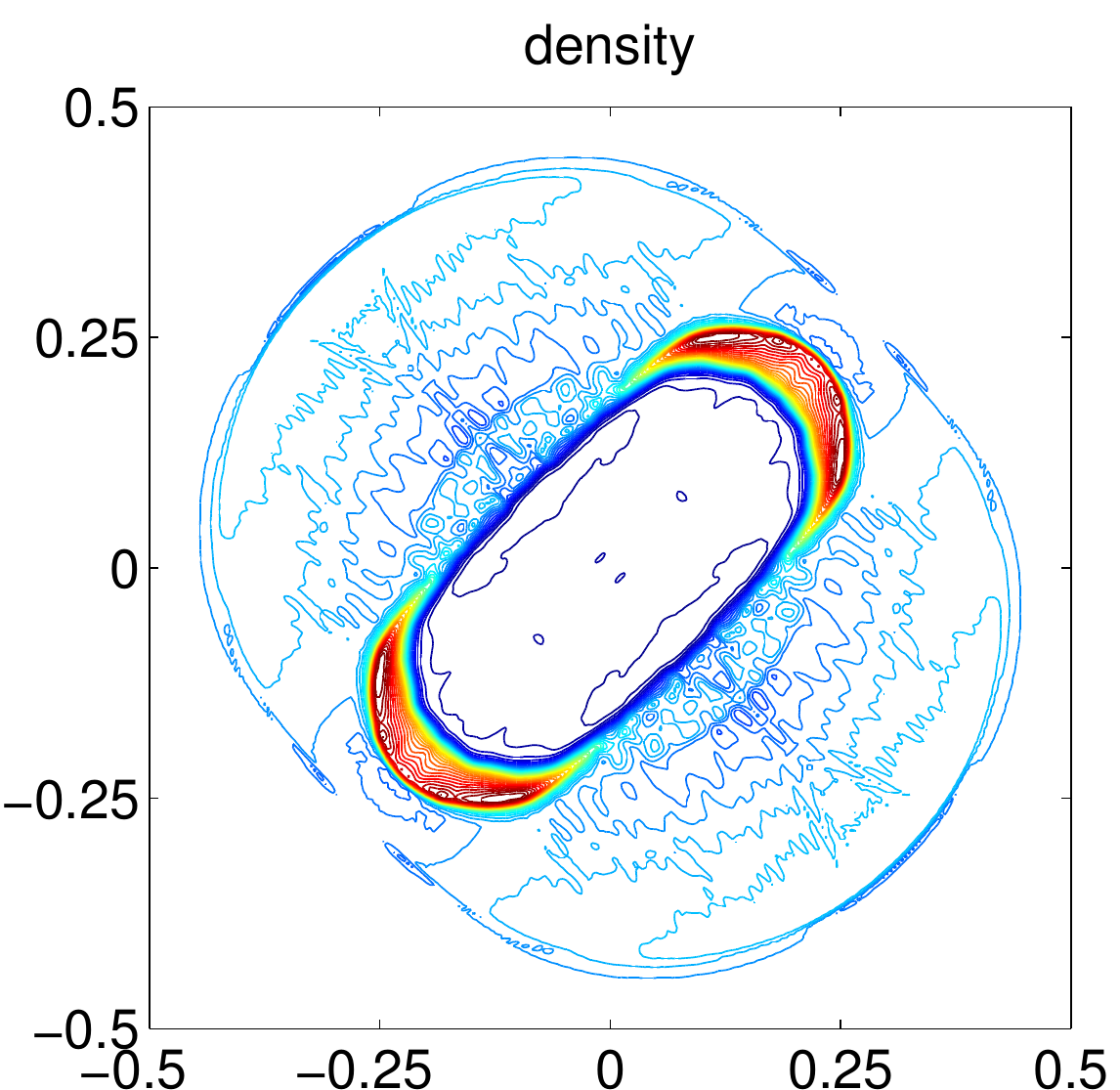} & 
	(b)\includegraphics[width=0.42\textwidth]{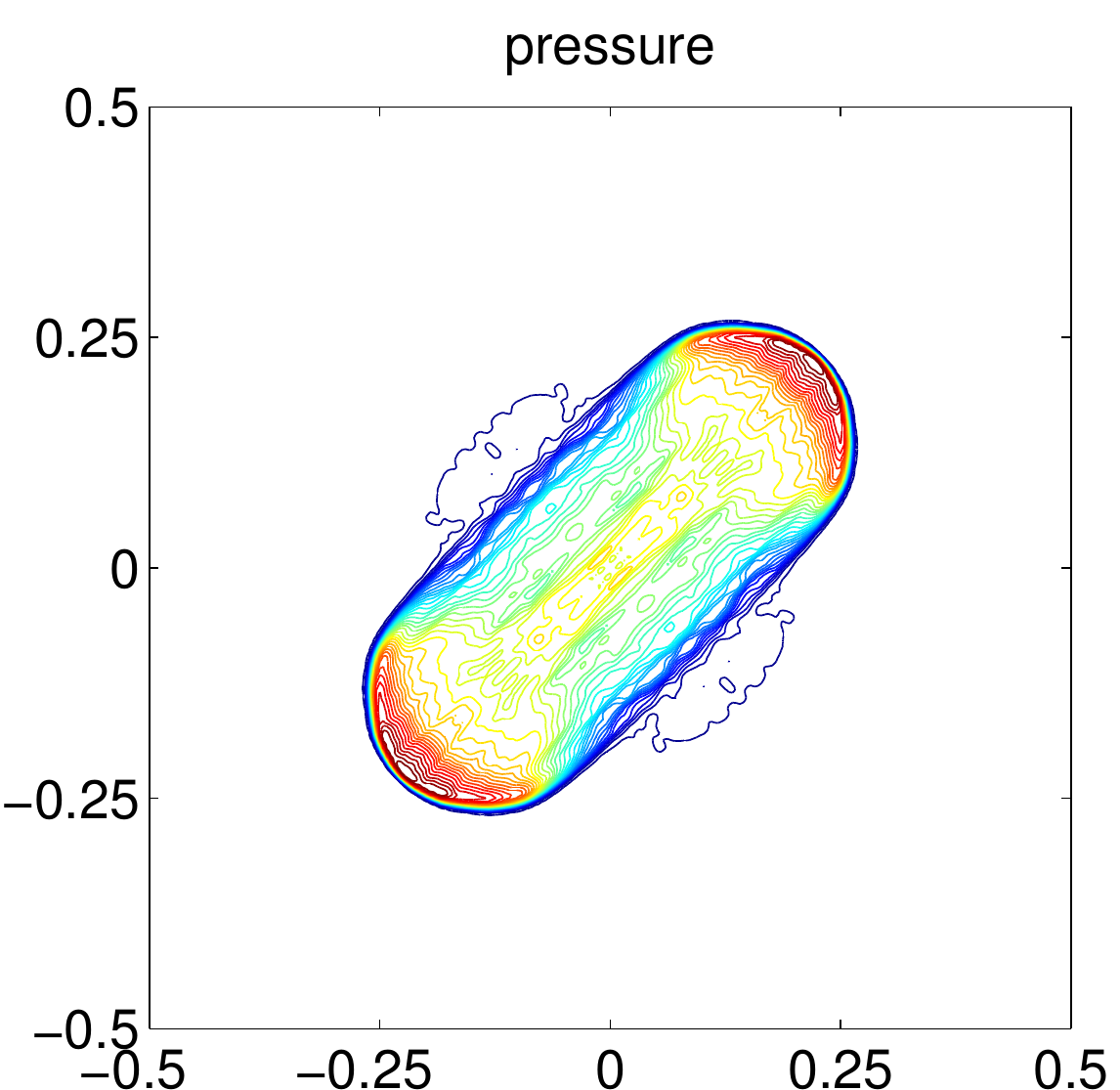} \\
	(c)\includegraphics[width=0.42\textwidth]{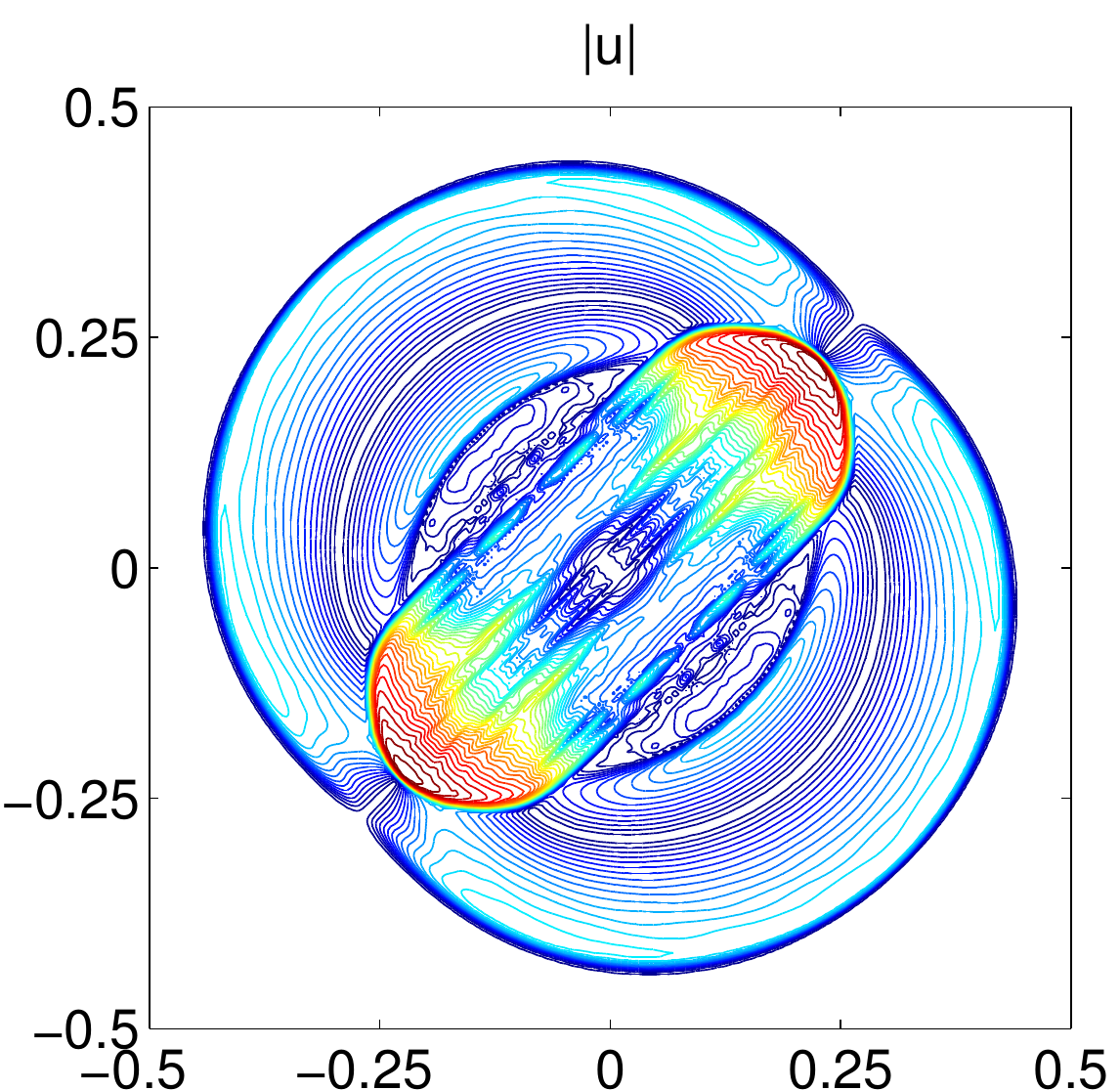} & 
	(d)\includegraphics[width=0.42\textwidth]{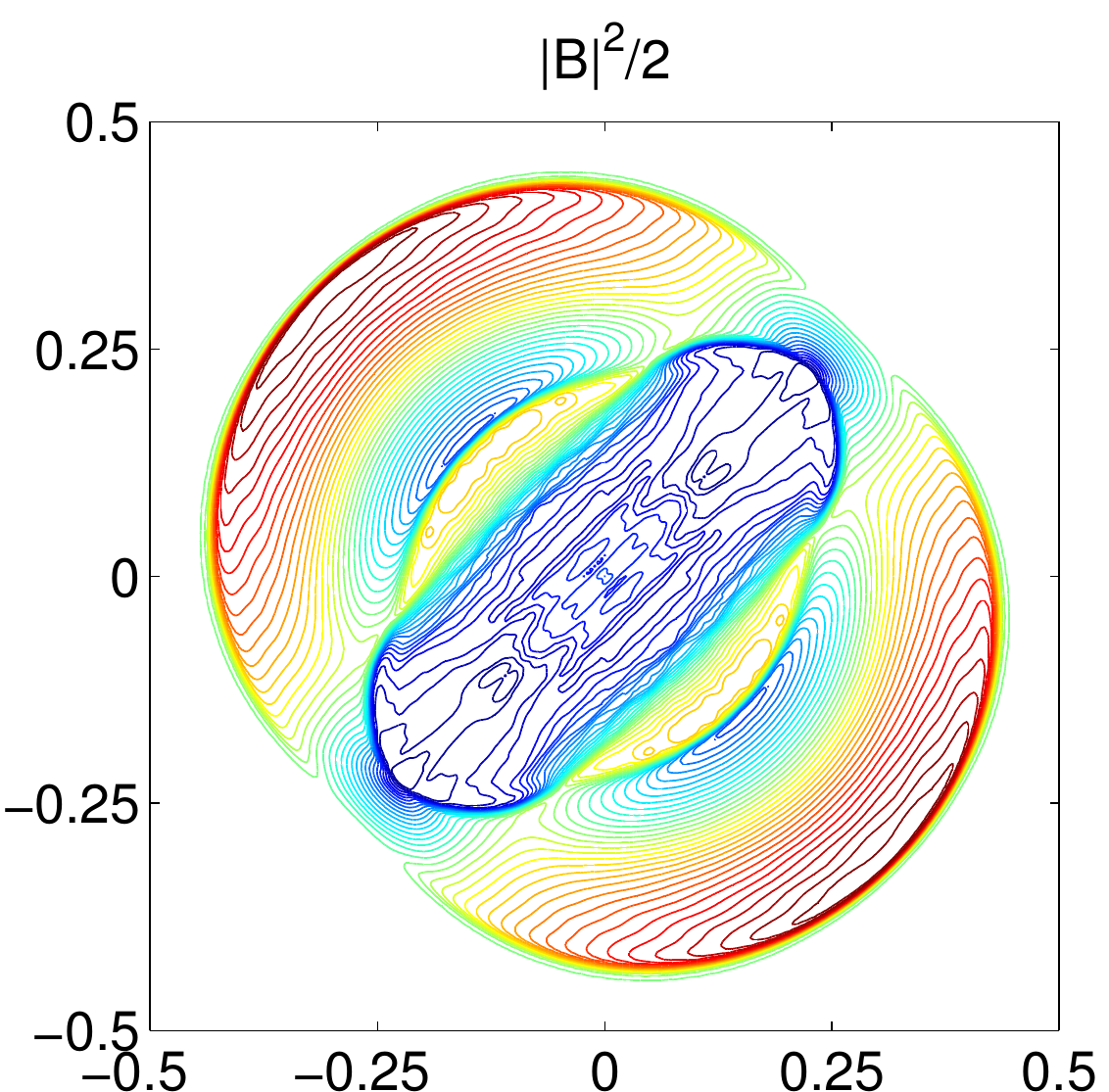}
\end{tabular}
      \caption{2D blast problem. 
      Shown in these panels are plots at time $t=0.01$ of 
      (a) the density,
      (b) the thermal pressure,
      (c) the norm of velocity and
      (d) the magnetic pressure. 
      40 equally spaced contours are used for each plot.
      The solution was obtained on a $256 \times 256$ mesh by WENO-HCL scheme.    
      \label{blast_hcl}}
\end{center}
\end{figure}

Another concern in the CT framework is that the energy is not conservative in our positivity-preserving WENO-CT scheme due to {\bf Option 2}. 
We also use this example to study this issue. 
We compare the results by the base WENO-HCL scheme and the 
positivity-preserving WENO-HCL scheme.
In Figure~\ref{blast_hcl}, we show the results by the base WENO-HCL scheme with the same resolution as Figure~\ref{blast2}.
The results look similar to those by the positivity-preserving WENO-CT scheme,
except there are some unphysical oscillations around the center region in Figure~\ref{blast_hcl}.
That is due to the divergence error in the base scheme. 
If we plot the divergence error in the time domain, 
we can clearly see the divergence error of the positivity-preserving WENO-CT scheme
stays around $10^{-12}$, while the error of the WENO-HCL scheme is around $10^0$ during the simulation.
Here the divergence error is defined as $L_1$-norm of $\divg \B$,
where the numerical $\divg$ operator is defined as a regular fourth-order central finite difference discretization.
As a common drawback in CT framework when {\bf Option 2} is chosen,
the correction step leads to a loss of the conservation of the total energy.
The results are plotted in Figure~\ref{blast_hcl}, where we can see 
the relative total energy error is around $10^{-3}$ 
while the conservative WENO-HCL scheme has an error about $10^{-12}$.
But we find this loss of total energy will decrease as the mesh is refined.
However, we remark that the conservation of the total energy is important for some problems, 
such as those involving nonlinear strong discontinuities.
A high-order positivity-preserving conservative scheme with  
the divergence error controlled will be part of our future work.
However, it is very difficulty, if not impossible, in the CT framework to satisfy 
all the requirements simultaneously.
A better way to control the divergence error is needed for this purpose.

\begin{table}
\begin{center}
\begin{Large}
\caption{Comparisons of different schemes solving the 2D blast problem. 
The column of ``Positivity'' lists if a negative solution is observed in the simulation. The column of ``Stability'' lists if the simulation run stably to $t = 0.01$.
In order to make a fair comparison, the positivity of density and pressure is only checked at each time step $t^n$. 
\label{tab:comp}}
\end{Large}
\begin{tabular}{|c || c c || c c|| c c || c c |}
  \hline
{\normalsize{ Mesh}} & \multicolumn{2}{c||}{WENO-HCL} & \multicolumn{2}{c||}{WENO-CT-OP1} & \multicolumn{2}{c||}{WENO-CT-OP2} & \multicolumn{2}{c|}{PP-WENO-CT-OP2} \\ 
\cline{2-5} \cline{6-9}
 & {\normalsize  Positivity}  &  {\normalsize  Stability} & {\normalsize Positivity } & {\normalsize Stability }
 & {\normalsize  Positivity}  &  {\normalsize  Stability} & {\normalsize Positivity } & {\normalsize Stability }
 \\ \hline \hline
{\normalsize $150 \times 150$}  & No	&	No	&	No	&	No
								& No	&	Yes	&	Yes	&	Yes\\
								
{\normalsize $200 \times 200$}  & No	&	Yes	&	No	& 	No
								& No	&	Yes	&	Yes & 	Yes\\
								
{\normalsize $256 \times 256$}  &No	&	Yes	&	No	&	No	
								&No	&	No	&	Yes &	Yes\\
 \hline
 \end{tabular}
 \end{center}
 \end{table}

\subsubsection{3D blast problem}
\begin{figure}
\begin{center}
\begin{tabular}{cc}
  	(a)\includegraphics[width=0.42\textwidth]{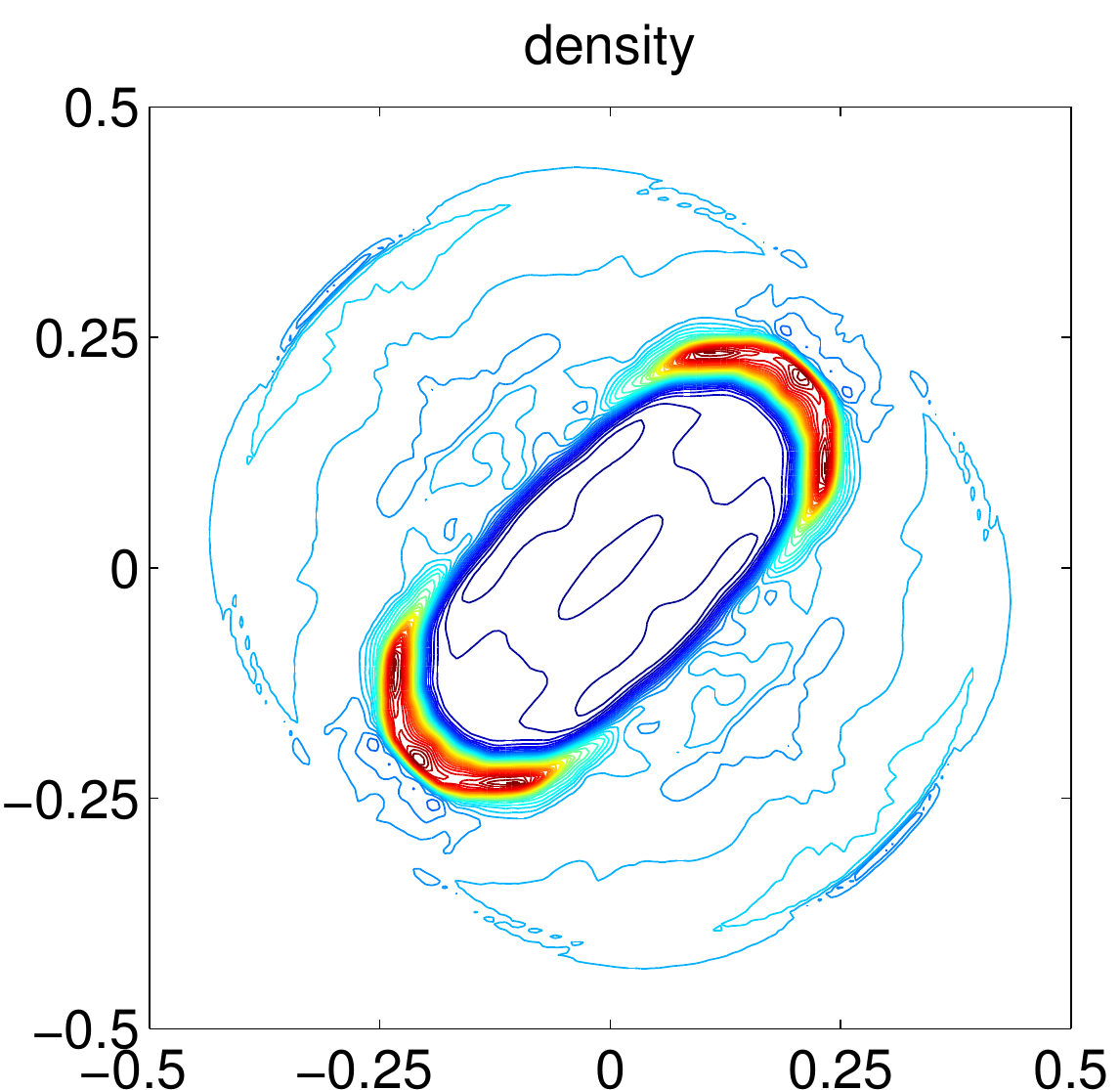} & 
	(b)\includegraphics[width=0.42\textwidth]{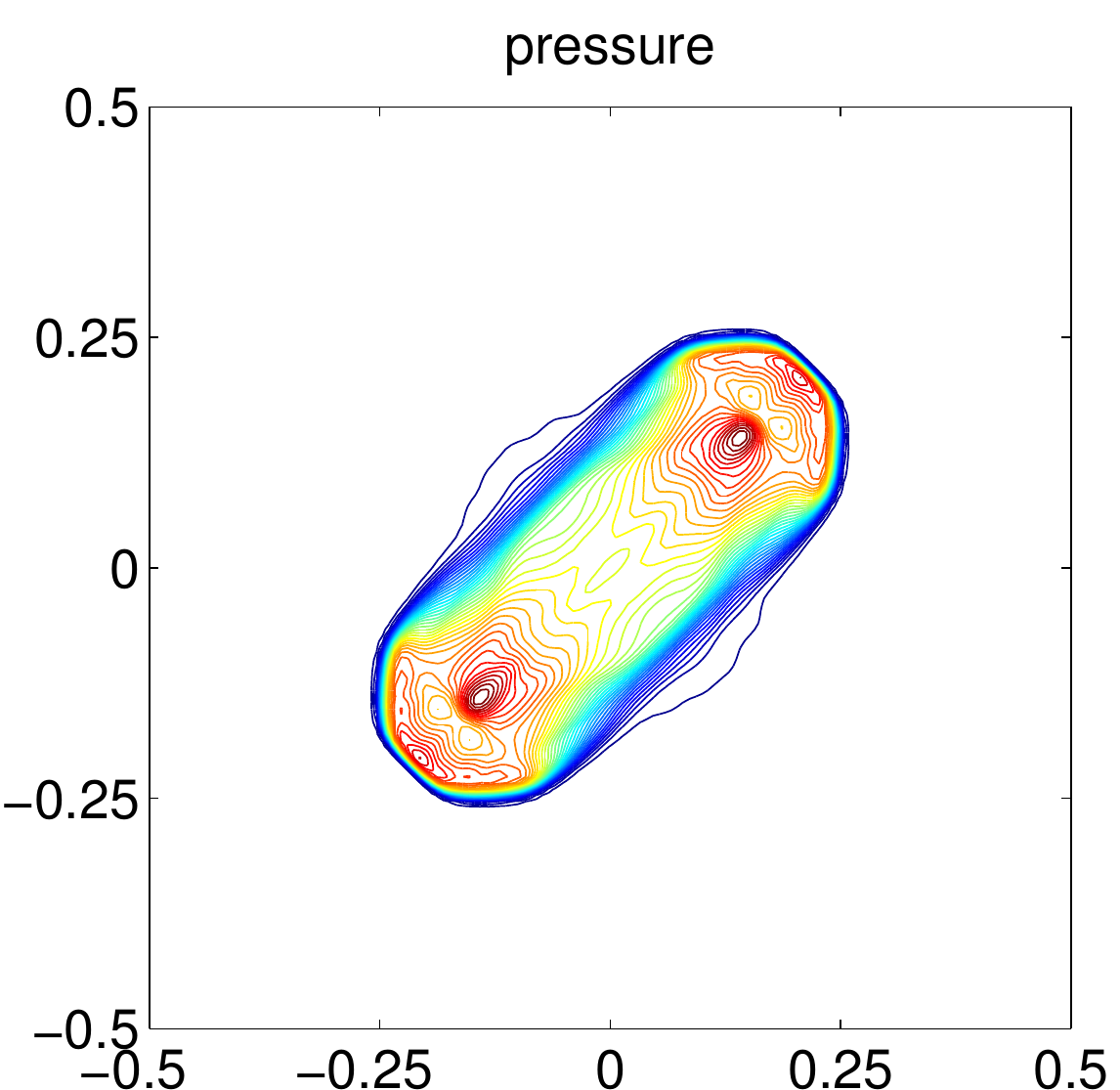} \\
	(c)\includegraphics[width=0.42\textwidth]{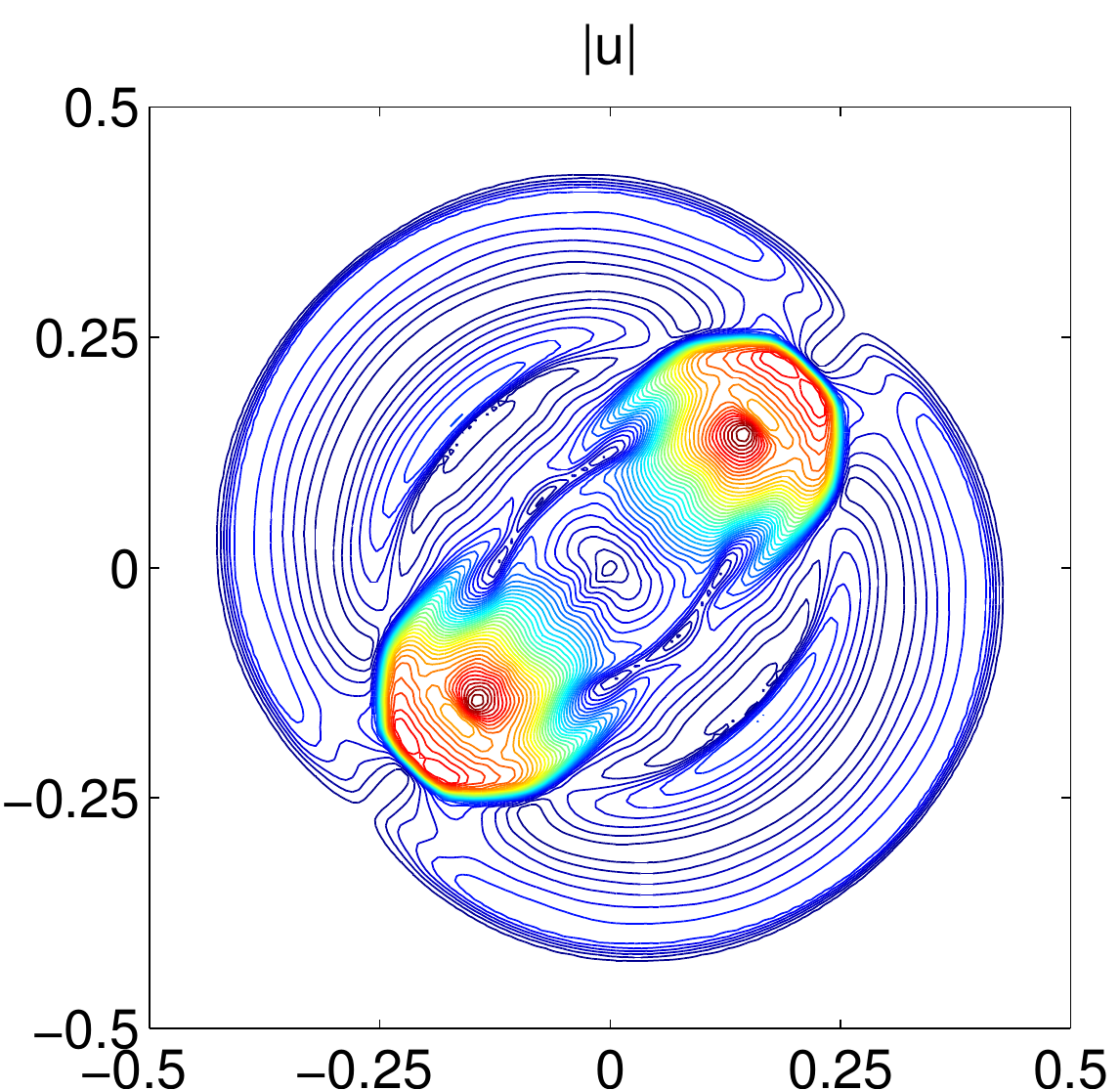} & 
	(d)\includegraphics[width=0.42\textwidth]{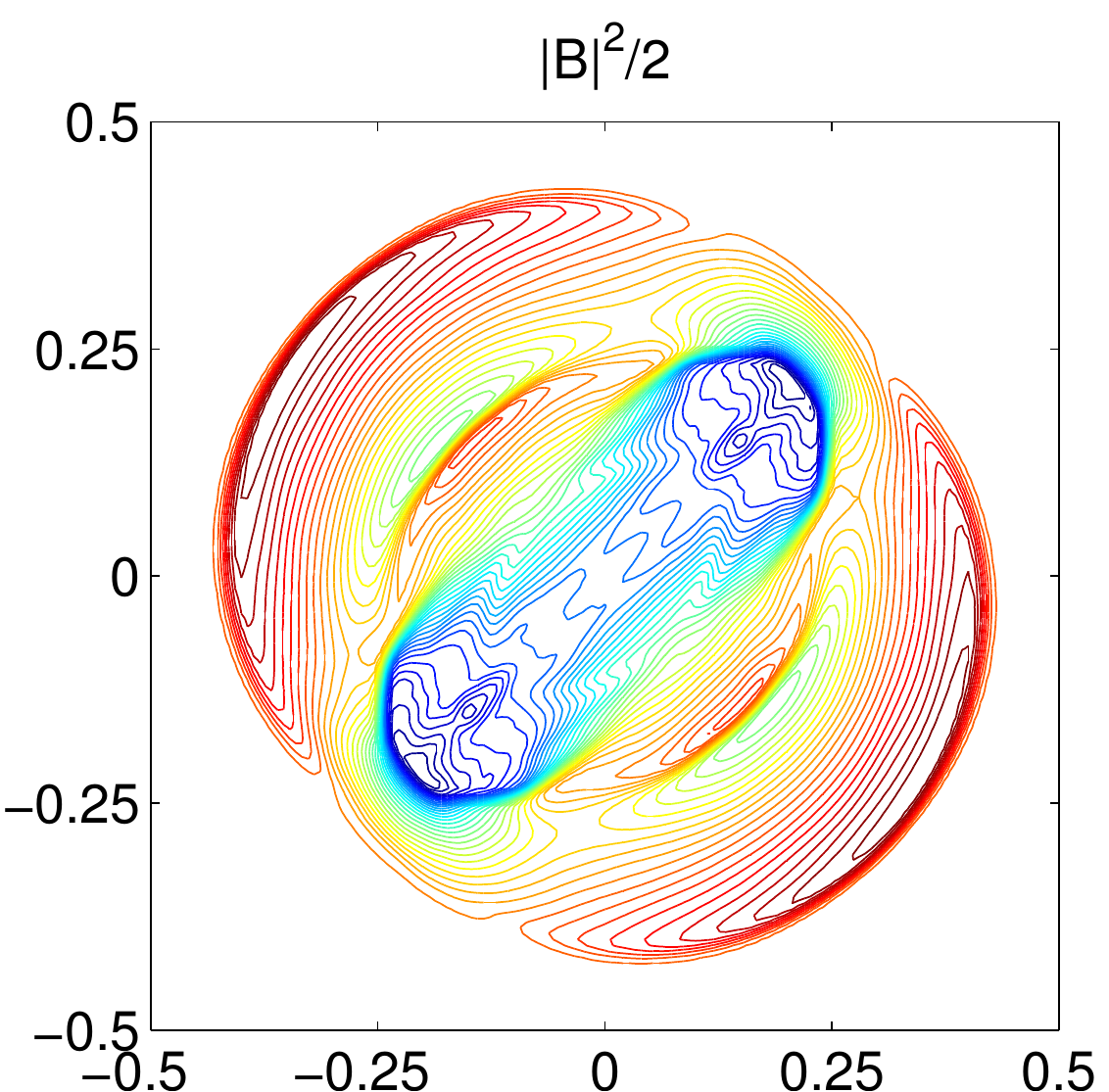}
\end{tabular}
      \caption{3D blast problem. 
      Shown in these panels are plots at time $t=0.01$ and cut at $z = 0$ of 
      (a) the density,
      (b) the thermal pressure,
      (c) the norm of velocity and
      (d) the magnetic pressure. 
      40 equally spaced contours are used for each plot.
      The solution was obtained on a $150 \times 150 \times 150$ mesh.     
      \label{blast3d}}
\end{center}
\end{figure}

\begin{figure}
\begin{center}
\begin{tabular}{cc}
(a) \includegraphics[width=0.42\textwidth]{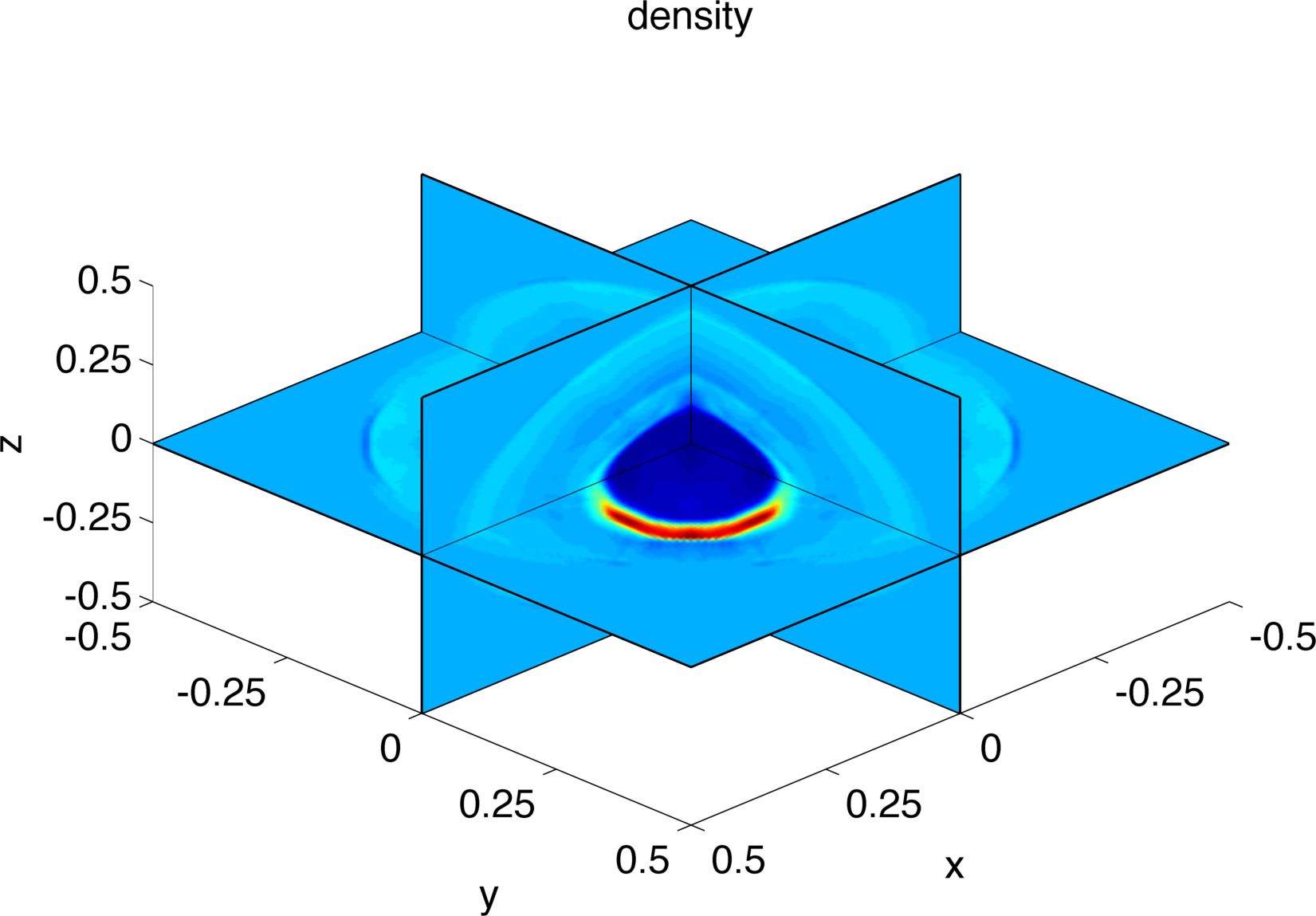} &
(b) \includegraphics[width=0.42\textwidth]{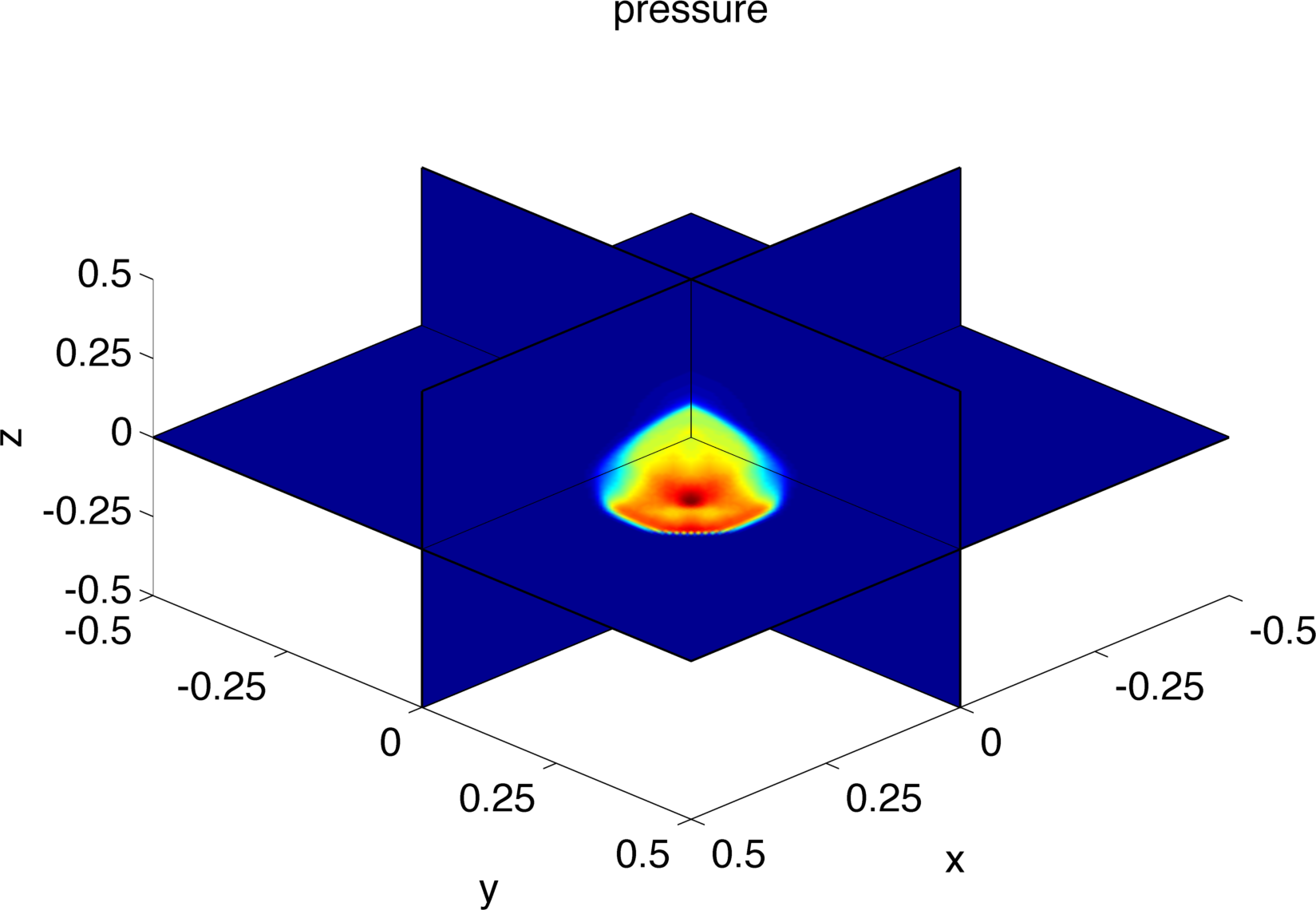}
\end{tabular}
\caption{3D blast problem. 
      Shown in these panels are 3D plots at time $t=0.01$ of 
      (a) the density and
      (b) the thermal pressure.
      The solution was obtained on a $150 \times 150 \times 150$ mesh.     
      \label{3dplot}
}
\end{center}
\end{figure}

The last problem we investigate is a fully 3D version of the blast problem.
It is used to test the behavior of the positivity-preserving WENO-CT3D scheme.
The initial conditions consist of an initial background:
\begin{equation}
(\rho, u_x, u_y, u_z, p, B_x, B_y, B_z) = (1, 0, 0, 0, 0.1, 100/\sqrt{4 \pi}/\sqrt{2},100/\sqrt{4 \pi}/\sqrt{2},0)
\end{equation}
and a spherical pressure pulse $p = 1000$ within a radius $r=0.1$ from the centre of the domain.
The initial conditions for the magnetic potential are
\begin{align}
\Av(0,x,y,z) = (0, 0, 100/\sqrt{4 \pi}/\sqrt{2} y - 100/\sqrt{4 \pi}/\sqrt{2} x ).
\end{align}
The computational domain is $[-0.5,0.5]^3$. 
Outflow boundary conditions are used on all the sides. The numerical simulation is performed on a $150 \times 150 \times 150$ mesh. 
In Figure \ref{blast3d} we show the results of the solutions cut at $z = 0$.
To distinguish this 3D case from the 2D blast case, 
we also present the 3D plots of the density and pressure in Figure \ref{3dplot}, 
which clearly indicates its spherical structures.
The solution is comparable to the 3D results in \cite{article:GaSt08, article:MiTz10, article:Zi04}.
We note that negative pressure is observed at time $t = 0.0033$ if the positivity preserving limiter is not applied.

\begin{figure}
\begin{center}
\begin{tabular}{cc}
  	(a)\includegraphics[width=0.42\textwidth]{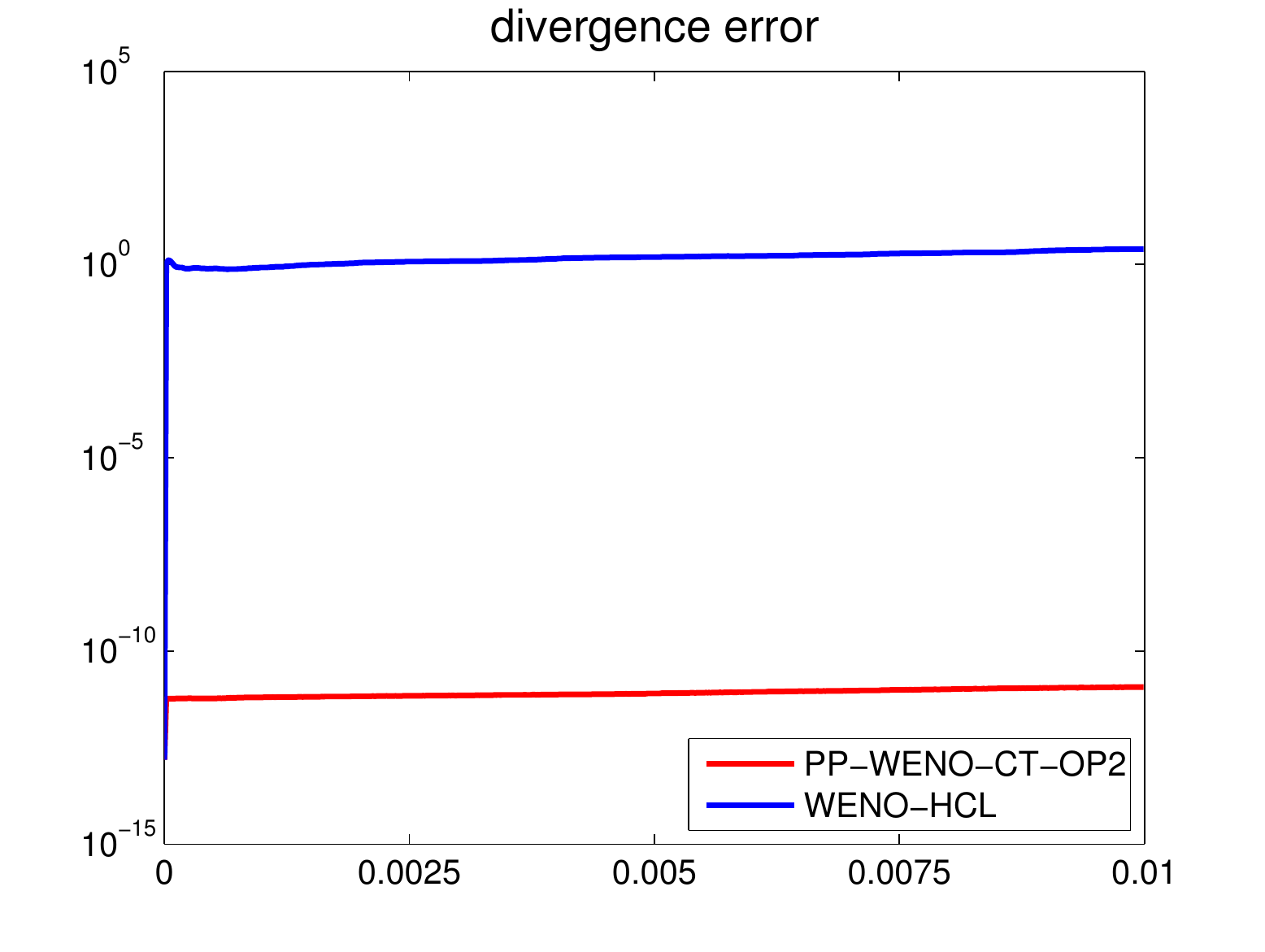} & 
	(b)\includegraphics[width=0.42\textwidth]{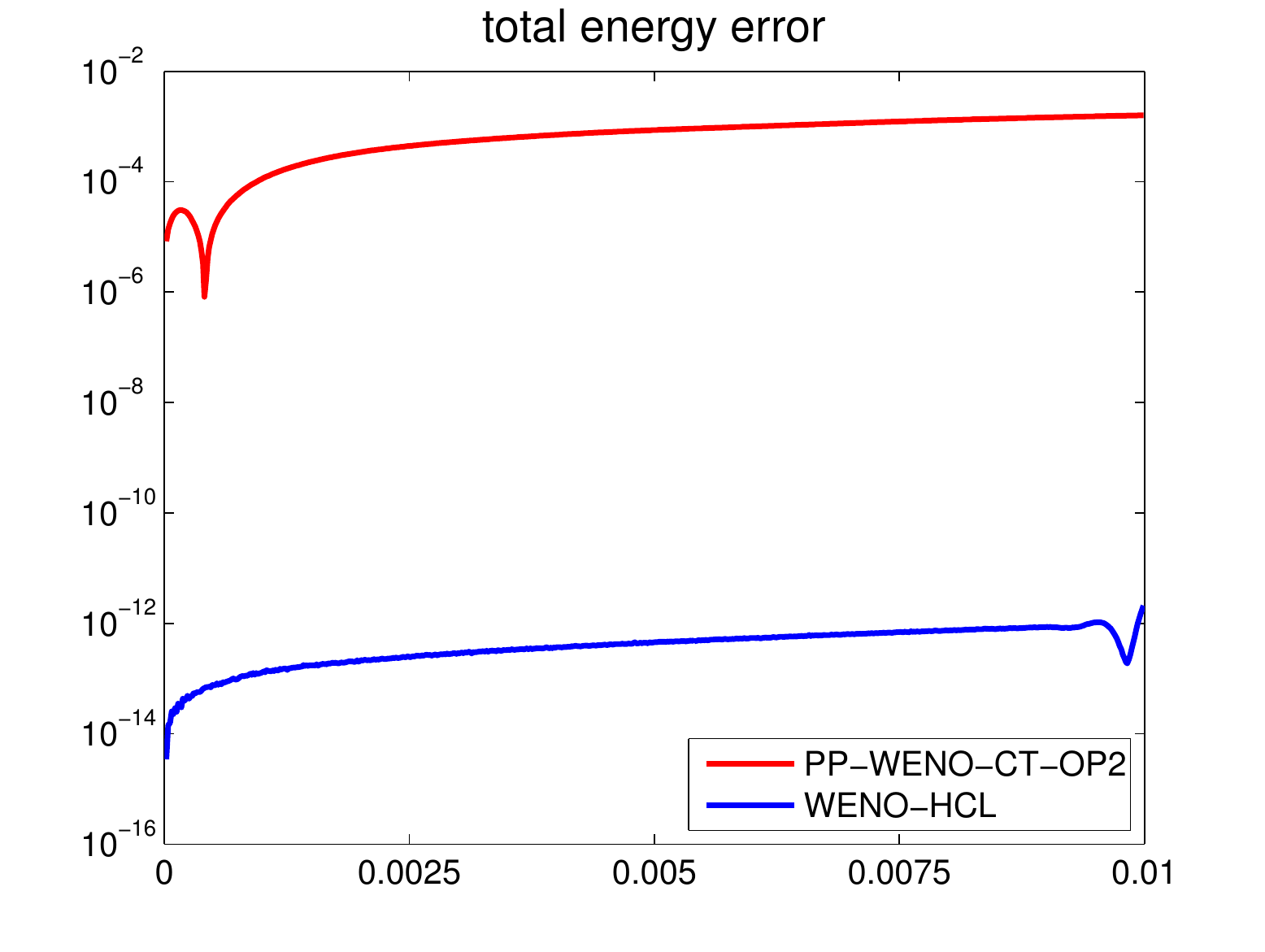}
\end{tabular}
      \caption{Comparisons between WENO-HCL scheme and positivity-preserving WENO-CT scheme for the 2D blast problem. 
      Shown in these panels are plots of 
      (a) the divergence error and 
      (b) the relative total energy error in the time domain $t \in [0,0.01]$.     
      \label{blast_comp}}
\end{center}
\end{figure}

\section{Conclusion}
\label{conclu}
In this paper we proposed a class of novel high-order positivity-preserving 
finite difference schemes for the 1D and multi-D ideal MHD systems. 
In the 1D case, a positivity-preserving limiting technique was applied to modify
high-order WENO-HCL flux with the first-order Lax-Fridrichs flux to produce positive density and pressure.
In multi-D cases, the limiting technique was also applied to the 
hyperbolic solver, followed by a constrained transport technique evolving the magnetic 
potential to control the divergence errors.
The main advantage of the proposed schemes is, the high order of accuracy, a discrete divergence-free condition and positivity of solutions can be attained  at the same time.
The overall scheme shares the same CFL constraint as the low-order Lax-Fridrichs scheme, 
without extra restriction resulting from the limiting process.
We demonstrated the effectiveness and efficiency of the positivity-preserving schemes by 1D, 2D and 3D numerical examples.
A strict proof for high-order accuracy of the proposed limiting technique will be part of our future work.

\bigskip

\noindent
{\bf Acknowledgements.}
AJC is supported by AFOSR grants FA9550-11-1-0281, FA9550-12-1-0343 and
FA9550-12-1-0455, NSF grant DMS-1115709, MSU Foundation grant SPG-RG100059, and by ORNL under an HPC LDRD. ZX is supported by NSF grant DMS-1316662.

\bibliographystyle{plain}
\bibliography{BigBib}

\end{document}